\documentclass{article}

\usepackage[english]{babel}

\usepackage[letterpaper,top=2cm,bottom=2cm,left=3cm,right=3cm,marginparwidth=1.75cm]{geometry}

\usepackage{amsmath}
\usepackage{graphicx}
\usepackage[colorlinks=true, allcolors=blue]{hyperref}

\title{A Driving Risk Surrogate and Its Application in Car-Following Scenario at Expressway}
\author{Renfei Wu,Linheng Li,Haotian Shi,Yikang Rui,Dong Ngoduy,Bin Ran}
\begin{document}
\maketitle

\begin{abstract}
Traffic safety is important in reducing death and building a harmonious society. In addition to studies of accident incidences, the perception of driving risk is significant in guiding the implementation of appropriate driving countermeasures. Risk assessment can be conducted in real-time for traffic safety due to the rapid development of communication technology and computing capabilities. This paper aims at the problems of difficult calibration and inconsistent thresholds in the existing risk assessment methods. It proposes a risk assessment model based on the potential field to quantify the driving risk of vehicles. Firstly, virtual energy is proposed as an attribute considering vehicle sizes and velocity. Secondly, the driving risk surrogate(DRS) is proposed based on potential field theory to describe the risk degree of vehicles. Risk factors are quantified by establishing submodels, including an interactive vehicle risk surrogate, a restrictions risk surrogate, and a speed risk surrogate. To unify the risk threshold, acceleration for implementation guidance is derived from the risk field strength. Finally, a naturalistic driving dataset in Nanjing, China, is selected, and 3063 pairs of following naturalistic trajectories are screened out. Based on that, the proposed model and other models use for comparisons are calibrated through the improved particle optimization algorithm. Simulations prove that the proposed model performs better than other algorithms in risk perception and response, car-following trajectory, and velocity estimation. In addition, the proposed model exhibits better car-following ability than existing car-following models.
\end{abstract}

\section{Introduction}
Road traffic accidents are among the top ten leading causes of death and morbidity globally \cite{blas2010equity}. Roughly 1.35 million people die in automobile accidents annually, with 20 to 50 million people suffering non-fatal injuries \cite{WHO}. Human errors are the primary cause of these traffic accidents \cite{nadimi2016calibration}. Risk assessment is of great significance for avoiding human errors, which has great potential to reduce traffic accidents. More advanced risk assessment approaches have recently been developed thanks to the fast development of communication technology, computing capabilities, and data gathering technologies. Specifically, real-time perception and data transmission between vehicles and intelligent facilities can be achieved through advanced communication technology, which enables information sharing regarding vehicles’ motion data (e.g., position, speed, and acceleration). Based on the real-time updated information, the driving risk of vehicles can be quantified in real-time to guide the operation of appropriate countermeasures.

Many researchers focused on risk assessment in recent years. Generally, risk assessment included two types of typical methods, the post-crash stage and the pre-crash stage. In the post-crash stage, risk assessments of road traffic were generally measured on the basis of statistical accident data \cite{arun_systematic_2021, schlogl2019methodological,liu_multivariate_2018,staton2016road}. Knowledge of the variables influencing the incidence of highway crashes can be expanded by taking into account spatial and temporal correlations, according to the review \cite{mannering_analytic_2014}. For example, Hojjati-Emami \cite{hojjati2013stochastic} presented the causes of vehicle accidents on roads using a fault tree model. Meanwhile, road safety indexes, the incidence of highway crashes, and other relative indexes were proposed to evaluate risk levels. They were usually obtained from statistical models such as accident regression analyses model \cite{chen2017monitoring,lord2010statistical}, extreme value theory (EVT)\cite{alozi_evaluating_2022}. 

 However, statistical analysis of collision data focused on crash events in the post-crash stage. It left out the majority of safe interactions that were common for naturalistic driving and intended to avoid collisions \cite{mohammadian_integrating_2021}. Furthermore, slight collisions were generally not included in traffic accident records, and traffic risk analysis methods based on previous accident data had relatively low adaptability and accuracy limitations. It was not enough to only research the crash incidents to guide driving safety. To prevent accidents, learning the risk threshold of safe interactions was an indispensable task. Pre-crash research intended to figure out the methodology of safe interactions from naturalistic driving data. In the pre-crash stage, classification was proposed on two aspects: physics-based models and potential field-based theories—which will be introduced as follows. Physics-based models were risk assessment methods based on time or kinematics. They were most commonly used in longitudinal analyses, such as the time-to-collision (TTC), the inverse of time-to-collision (TTCi), deceleration rate to avoid a crash (DRAC) \cite{qu2014safety,ward2015extending,sharifi2016time}, and lateral analyses, such as the vehicle’s current position (CCP), time-to-lane-cross (TLC), and the variable rumble strip (VRBS) \cite{tan2017use,sun2011elevated,wang2010applying}. However, in real traffic environments, traffic risks were continuous variables that cannot be artificially classified into longitudinal or lateral directions. Li \cite{li2020threat} and Lu \cite{lu2021performance} compared the performances of different kind of pre-crash risk assessment measures and their discussions showed that the physics-based modeled implement one-dimensional risk assessment and performed poorly in timeliness and robustness, while potential field-based models performed better in robustness and demonstrated a two-dimensional risk assessment but were computationally demanding.

Potential field-based theory for road traffic risk assessments had been widely used in obstacle avoidance and path planning for intelligent vehicles in recent years. They could represent risks in two or three dimensions because of their advantages in terms of environmental description. Some researchers introduced gravity and repulsion into potential field theory. Krogh \cite{krogh1986integrated} firstly applied the potential field theory to the perception of the surrounding environment of a robot. Gravity was used to indicate the walkable area, and repulsive force indicated the non-walking area. This was mainly used for the path planning decisions of a robot. Hsu \cite{hsu2012conceptual} introduced the concept of the gravitational field to establish a car-following model based on Krogh's research. Ni \cite{ni2013unified} proposed the field theory of a microscopic traffic model, which was mainly applied to the car-following model. He related the vehicle safety forward to a combination of a free fall and repulsive forces brought by surrounding vehicles. Some researchers applied the electronic field and spring to present the driving risk. Reichardt \cite{reichardt1994collision} built a vehicle-centered power electronic field combined with the surrounding environment. However, although the theory was in line with the general scenarios, there were still many realistic scenarios that could not be matched. To solve the problem, Tsourveloudis \cite{tsourveloudis2001autonomous} combined the electrostatic potential field theory with two-layer fuzzy logic reasoning. The environmental road network was analogized to a resistance network, in which each vehicle generated an electrostatic potential field. The optimal path of the vehicle was planned by the current in the network, thereby realizing the automatic path planning of an intelligent vehicle in a dynamic environment. Sattel \cite{sattel2008robotics} combined potential field with elastic order from the field of robotics. The elastic order was used to calculate the minimum path in the potential field to plan a vehicle path and avoid falling into the local optimum. Matsumi \cite{matsumi2013autonomous} examined an autonomous collision avoidance system based on electric vehicle braking torque. In the system, the braking maneuver intensity was determined using potential field theory, which considered the potential hazards because of occlusions in the intersection. Raksincharoensak \cite{raksincharoensak2013predictive}proposed a virtual spring model that connected a vehicle and pedestrian to determine the potential field in unsignaled intersections with poor visibility. The traffic risk was evaluated based on simple physical characteristics, while it involved complex factors. Thus, Wang \cite{wang2015driving} constructed a unified Driving Safety Field (DSF) model with an elliptic structure, and the DSF includes a potential field related to road lines, a kinetic field about surrounding vehicles, and a behavior field that was determined by the individual driver characteristics. Rasekhipour \cite{rasekhipour2016potential} assigned different field functions to different types of obstacles and road structures for path planning. Moreover, vehicle collision avoidance was in combination with the planning path and vehicle dynamics constraints. Zheng \cite{zheng2018novel} developed a model based on equivalent force, which calculated the road traffic risk field in different directions. Huang \cite{huang2020probabilistic} proposed a probabilistic risk assessment framework combined with intention identification to detect a threat in advance and apply it to lane-changing scenarios. Li \cite{li2020dynamic} researched a dynamic risk surrogate and conducted its performance in car-following scenarios. In all these studies, the road traffic risk was quantified and presented as a distribution in a two-dimensional risk map, in which the potential energy value denoted the road traffic risk level. 

Potential field-based theory for road traffic risk assessments intended to serve vehicle decisions. Usually, vehicles made micro-driving decisions according to assessment results of the risk concerns of their surroundings. For example, if there was a potential collision in the car-following scenario, the following vehicle would normally change lanes or decelerate as the leading vehicle decelerated to keep safe. To put it another way, risk assessment directly impacted microscopic vehicle behavior, and the driving risk surrogate reflected this risk status in the form of field strength. Despite distinct risk models' diverse structures and expressions due to different modeling views, there were no clear differences in assessment mechanisms. These models were based on actual vehicle operation to accurately define microscopic behavior and traffic flow features through model analysis. The main distinction lay in the implemented approaches supported by different transportation theories. The potential field theories integrated varied traffic factors, which led to a challenging calibration process. Among these factors, vehicles' masses were normally designed as different constants in existing potential field theories. Thus, the virtual mass used to be calibrated offline and processed ideally. However, vehicles were heterogeneous in naturalistic driving environments. The masses of vehicles were different from each other even for the same type of vehicles. Furthermore, the masses of vehicles could not be exhaustive precisely. The ubiquitous cameras and radars roadside hardly perceived the mass, although they could collect a lot of vehicle status information. Moreover, the risk thresholds of field energy were not unified in different scenarios. The thresholds were fixed in the related interactions scenario. The threshold needed to change when the interaction was transferred to another type. That meant effective risk assessment depended on enumeration and identification of interaction scenarios. However, there were tens of thousands of interactive scenes. And it was difficult to recognize scenario type accurately in a very short time.

To fill the research gaps, this paper develops an innovative driving risk surrogate(DRS) to solve challenging calibrations and inconsistent threshold problems. In the modeling process of DRS, we consider sizes as the critical vehicle attribute, while that was mass in the previous models. Vehicle sizes are more intuitive and more accessible to obtain than mass through real-time perception. The specific size information enables the parameters of the virtual energy of heterogeneous vehicles to be trained and calibrated by the real driving dataset. A relationship between field energy and acceleration is proposed for the problem of inconsistent thresholds. The acceleration is used to calibrate parameters in different scenarios according to naturalistic driving data. Also, we propose a speed risk surrogate to demonstrate the driver's pursuit of efficiency.

The proposed DRS accurately characterizes the driving risk level and explains varied microscopic driving behaviors. The value of DRS directly impacts vehicle behavior at the microscopic level. Since driving safety circumstances affect how a driver executes their next action, microscopic behaviors can be guided by the risk strength of DRS. Based on the model establishment, DRS is applied to the car-following model. Results show that DRS can successfully present the rules of microscopic car-following behavior and perform well in trajectory estimation.

The main contributions of this paper to the literature are threefold:

1. The driving risk assessment issue is tackled from a new perspective: an integrated model is developed by interactive vehicle risk, restrictions risk, and speed risk. 

2. The attribution of vehicles, virtual energy, is introduced based on size and velocity. It makes DRS more in line with the risk varying in drivers' perception and solves the problem of idealized processing mass in conventional models. 

3. Acceleration is designed as the output of DRS, solving the problem of difficult calibration and inconsistent thresholds. It is derived to motion status, and DRS is applied to simulate car-following behavior successfully.

The rest of the paper is organized as follows. Section 2 establishes DRS, and acceleration is deduced from the risk field strength to prepare for trajectory estimation. Section 3 introduces the naturalistic driving dataset, data pre-processing, and parameter calibration. In section 4, we present the risk response performance of four different models, including the proposed model. Simulations of trajectory estimation are conducted, and the results are analyzed and discussed. Finally, we conclude the paper in Section 5 with recommended future research.

\section{Driving Risk Surrogate Model}
The safety of a driving vehicle is influenced by its surroundings. Thus, it is necessary to analyze surroundings and classify them into different categories for risk modeling. In subsection 2.1, the factors are divided into three categories, including interactive vehicle risk, restriction objects risk, and desired speed risk. Vehicle attribute is designed in subsection 2.2 for representing important characteristics of each vehicle. For example, the risk formed by vehicles with different types is obviously different from each other even in the same motion state, driving scenario, and desired speed. In the following three subsections, the three types of risk surrogates are analyzed and proposed for the DRS modeling.

\subsection{Classification of Driving Risk Factors}
 Some studies indicated that the majority of road and traffic elements have a significant impact on driving behavior \cite{paschalidis_deriving_2020}. Human driver factors affecting driving risk can be divided into four categories: physiological-psychological factors, cognitive level, driving skills, and violations. Drivers will adopt a variety of driving behaviors to preserve their own safety. In general, a vehicle's goal, regardless of the scenarios, is to reach its desired speed while maintaining a safe clearance between itself and other vehicles.

Specifically, the factors affect driving safety, such as human driver, subject vehicle’s performance status, other road users, surrounding vehicles’ movement status, obstacles, traffic rules, and road alignment. Equipment failure, such as tire blowout, or brake failure, can affect vehicle performance status. This paper assumes that the performance status of the driver and the vehicle is at a normal level. Each driver with mental health and driving in well performance vehicle would try to avoid crashes and reach higher efficiency in expressways. In urban road scenes, vehicles during driving will interact with other road users, such as pedestrians and non-motor vehicles, and they will also impact driving safety. However, those road users do not appear on the expressway.

\begin{figure}[!t]\centering
    \includegraphics[width=16.0cm]{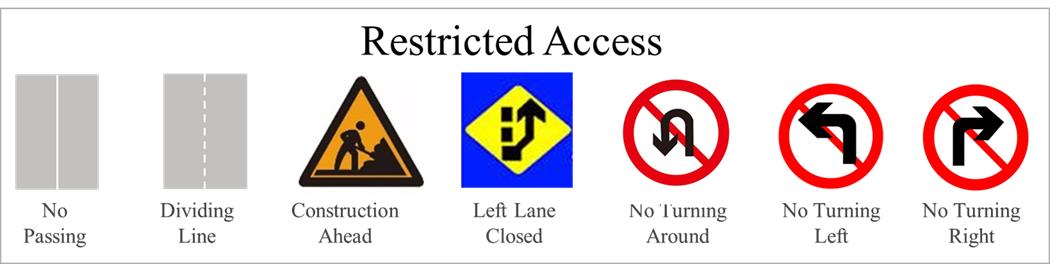}
    \caption{Categories of restricted access.}\label{Categories}
\end{figure}

As mentioned above, driving safety is also related to surrounding vehicles’ motion status, road environment, traffic rules, and road alignment. Undoubtedly, a driving vehicle interact and avoid conflict with other vehicles around it. And then, road environment refers to the static obstacles on or on both sides of the road. The drivers must avoid hitting these objects. Common objects of this kind include the isolation belt in the center of the road, the anti-collision pier at the edge of the road, the area of road construction, etc. Road markings and road signs define traffic rules for vehicle driving. The rich traffic rules related to driving safety are filtered out and are shown in \mbox{Fig. \ref{Categories}}. The road alignment also belongs to restricted access because of the same characteristic.

As a result, DRS can be separated into three categories: interactive vehicle risk surrogate, restrictions risk surrogate, and speed risk surrogate, which are described by $\overrightarrow {E_{ji}^v}$, $\overrightarrow {E_{wi}^t}$ and $\overrightarrow {E_i^s}$. They are field functions and act as field vectors from different risk sources. Due to varied directions of field vectors, there are phenomena in which the field vectors cancel others. However, for road safety, such phenomena are not consistent with reality. After superimposing diverse driving risk field vectors, the driving risk will normally grow rather than decrease. As a result, the total of the field strength modes created by distinct risk surrogates is used to describe the degree of the driving risk of a subject vehicle. Thus, the total risk field strength of DRS can be expressed as \mbox{Eq. (\ref{EQ_1})}.
\begin{equation}
\overrightarrow{E_{i}^{risk}}=\sum_{j} \overrightarrow{E_{ji}^{v}}+\sum_{w} \overrightarrow{E_{wi}^{t}}+\overrightarrow{E_{i}^{s}},
\label{EQ_1}
\end{equation}
where ${i}$ represents the subject vehicle, ${j}$ represents the interactive vehicle of the subject vehicle, ${w}$ is restriction of physical limitation or road rule, $ {E_i^{risk}}$ is the driving risk field strength and represents the sum of risks which ${i}$ suffered, $ {E_{ji}^v}$ is the interactive vehicle risk field strength which is the risk of ${i}$ suffered from interactive vehicle ${j}$, $ {E_{wi}^t}$ is the restrictions risk field strength and represents the risk of ${i}$ suffered from physical limitation or road rule ${w}$, and $ {E_i^s}$ is the speed risk field strength of ${i}$ in the current traffic flow. $\overrightarrow {E_i^{risk}}$, $\overrightarrow {E_{ji}^v}$, $\overrightarrow {E_{wi}^t}$, and $\overrightarrow {E_i^s}$ with over-arrows are field vectors.

\begin{figure}[!t]\centering
    \includegraphics[width=16.60cm]{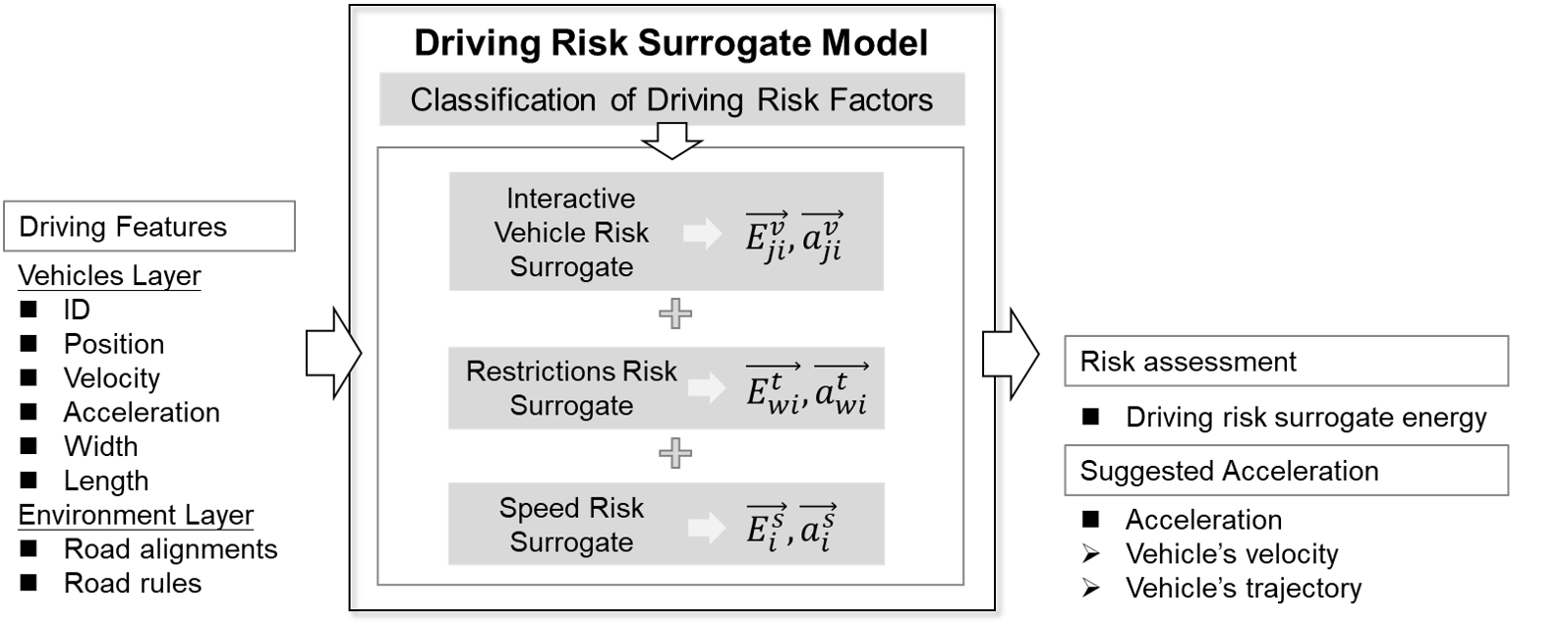}
    \caption{Block diagram of the proposed model.}\label{Block diagram}
\end{figure}

\mbox{Fig. \ref{Block diagram}} shows an overview of the DRS model. Driving features, including the vehicles layer and the environment layer, are input. The vehicle's layer comprises the ID, position, velocity, acceleration, width, and length of the subject vehicle and its neighboring vehicles. For privacy, vehicles are numbered with ID instead of license plates. The environment layer is composed of road layout and traffic rules. And then, the DRS model consists of the interactive vehicle risk surrogate, the restrictions risk surrogate, and the speed risk surrogate and outputs the result of the risk assessment and acceleration plan. The trajectories containing speed and location can already be inferred based on the acceleration.  
\subsection{Vehicle Attributes}

In previous research\cite{li2020risk}, the mass and type of a vehicle are the characteristics of the vehicle attributes. However, it does not match reality in some scenarios. For example, a truck with 200 tons presents a 200 times higher risk than a vehicle with two tons. This may not be intuitive enough. In particular, it is impractical that safety clearance between trucks with 200 tons requires 200 times as much space as automobiles with two tons. Drivers usually do not accurately perceive the mass of the interacting vehicle during driving. The most sensitive thing for drivers is the sizes of the interacting vehicles, including length, width, and height. Therefore, the sizes are selected to replace the mass and type and become an important attribute, which has been more in line with the actual effect. However, obtaining the vehicle's height information from the real driving dataset is difficult. According to vehicle types in the market, it can be found that a vehicle's height is related to its length and width (e.g., in general, the height of trucks is greater than that of SUVs (sport utility vehicles). The height of SUVs is usually greater than that of vehicles, which is similar to sizes of them). As a result, the vehicle's length and width are key attributes. 

In addition, Wang \cite{wang2016driving} suggested that the driving risk of automobiles changed at various speeds, even if they had the same mass. A combination of mass and velocity was described as equivalent mass. If a vehicle's speed increases, driving risk will increase synchronously in front and rear of the vehicle. However, as a driver, this is different from the risk perception when driving. With the increase of velocity, the risk in front of the vehicle increases, but that behind reduces. To overcome the identified issues of the current method, this paper proposes a new concept, namely "virtual energy," which is related to vehicle sizes, velocity, and velocity direction to reflect the vehicle's risk energy at a certain velocity. The mathematical expression of the virtual energy is defined as \mbox{Eq. (\ref{EQ_2})}.
\begin{equation}
Q_i=L_i^{s_l} W_i^{s_w} e^{\alpha v_i cos{\varphi} },
\label{EQ_2}
\end{equation}
where $Q_i$ is the virtual energy of vehicle $i$, $L_i$ and $W_i$ are the length and width of vehicle $i$, $v_i$ is the velocity, and $\varphi$ is the clockwise angle of any point to the velocity direction of the vehicle, $s_l$ means a power parameter of $L_i$, while $s_w$ is that of $W_i$, and $\alpha$ are undetermined parameters related to vehicle's velocity.

\subsection{Interactive Vehicle Risk Surrogate}
The interactive vehicle risk surrogate is derived from a moving vehicle. The strength of a subject vehicle's risk field is governed by its attributes, motion state, and distance to risk resources. In the interactive vehicle risk surrogate, the field strength at the position of the subject vehicle is related to the distance from the interactive vehicle. Let the space coordinates of the interactive vehicle's centroid be $(x_j,y_j)$. The distance $k_{ji}$ from any point $(x_{i},y_{i})$ to the vehicle's centroid in space can be expressed as \mbox{Eq. (\ref{EQ_3})}.
\begin{equation}
k_{ji}=\sqrt {(x_{i}-x_j)^2+(y_{i}-y_j )^2}.
\label{EQ_3}
\end{equation}

Suppose $k_{ji}$ is used to represent the negative correlation between distance and driving risk directly. In that case, the risk degree is the same at equidistant points, whether the vehicle is in the same lane as the target vehicle or another lane. However, this condition is inconsistent with reality. In practice, a vehicle should stay away behind a leading vehicle in its lane but is allowed relatively near the obstacle to the side or at the road boundary. For example, the degree of risk is similar when the vehicle is 30 meters behind the object vehicle and 2 meters to the side of the object vehicle. This is because there is no risk of collision between vehicles in adjacent lanes. Therefore, to better describe the relationship between the risk field strength and the distance, pseudo-distance is adopted in \cite{li2020novel}. The distance in real space is altered to describe the change of risk degree better when a subject vehicle approaches the object vehicle from different angles. 

\begin{figure}[!t]\centering
    \includegraphics[width=14.0cm]{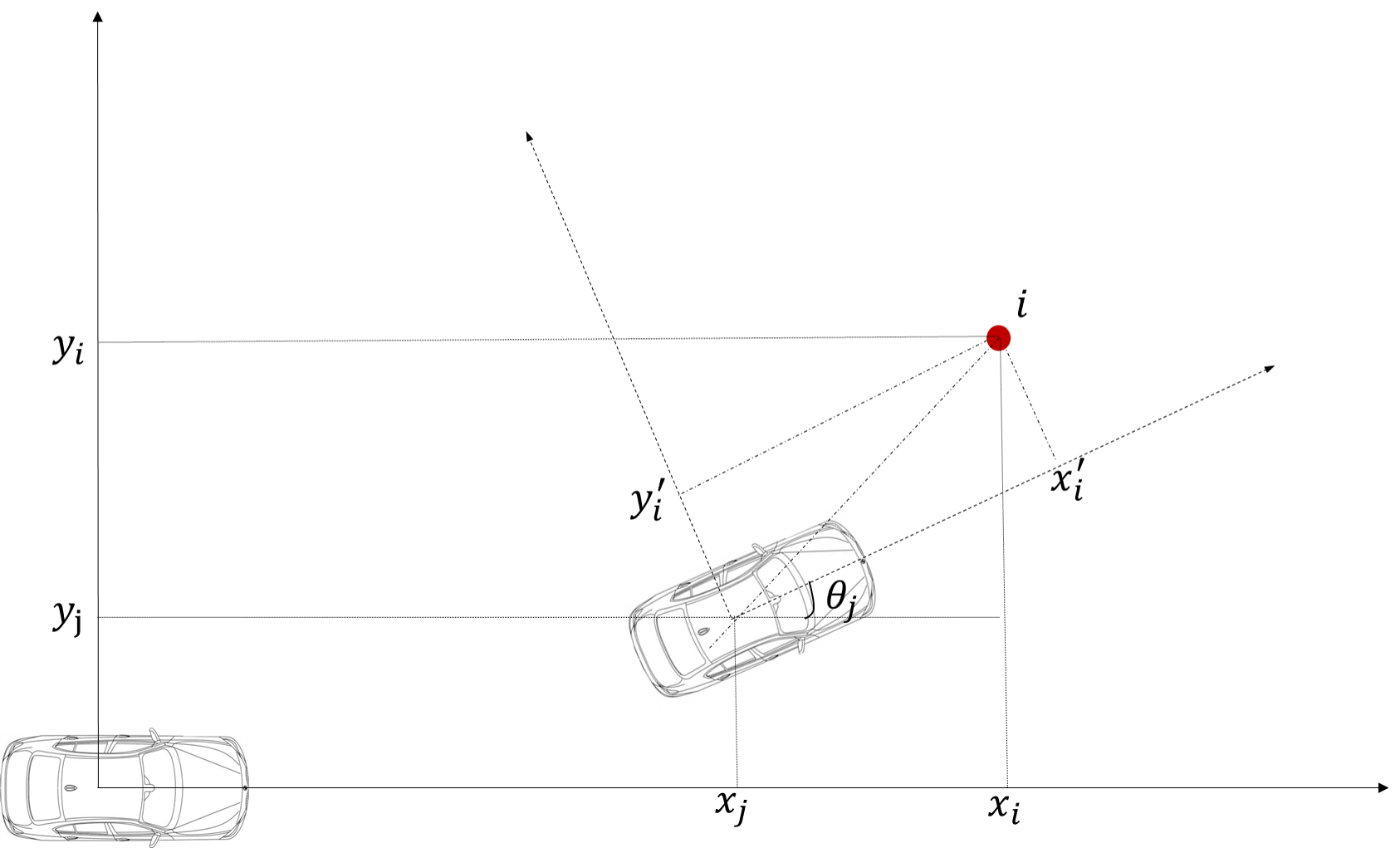}
    \caption{Sketch map of transforming global coordinate to vehicle coordinate.}\label{Sketch map}
\end{figure}

Based on the analysis above, virtual distance ${k_{ji}}^\prime $ can be expressed as \mbox{Eq. (\ref{EQ_4})}.
\begin{equation}
\begin{gathered}
\left\{\begin{array}{c}
\left.k_{ji}^{\prime}=\sqrt{\left[x_{i}^{\prime} \frac{e^{w_{1} v_{j}}}{1+w_{2} v_{j}}\right]^{2}+y_{i}^{\prime2} }\right.,\\
{\left[\begin{array}{l}
x_{i}^{\prime} \\
y_{i}^{\prime}
\end{array}\right]=\left[\begin{array}{cc}
\cos \theta_{j} & \sin \theta_{j} \\
-\sin \theta_{j} & \cos \theta_{j}
\end{array}\right]\left[\begin{array}{c}
x_{i}-x_{j} \\
y_{i}-y_{j}
\end{array}\right]},
\end{array}\right.
\label{EQ_4}
\end{gathered}
\end{equation}

where, $(x_{i}^{\prime},y_{i}^{\prime} )$ represents the value in the vehicle coordinate and $(x_{i},y_{i})$ is the original global coordinate, $(x_j,y_j)$ is the global coordinate of the object vehicle $j$, \mbox{Fig. \ref{Sketch map}} shows the sketch map transforming global coordinate to vehicle coordinate and $v_j$ is the velocity of the object vehicle, $\theta_j$ is the steering angle, $w_1$ and $w_2$ are undetermined parameters of velocity related to Doppler effect.

The interactive vehicle risk surrogate is a typical potential field applied in the traffic environment. The range of field for a vehicle in the driving process is restricted. The vehicle is only affected by a few vehicles surrounding it. With increasing distance, the impact will soon diminish. The phenomenon is similar to the characteristic of the short-range interaction between particles in physics. A potential field describes the interaction between two vehicles in the same way as it describes the interaction between two interacting particles. Inspired by this, we develop a vehicle risk surrogate that includes the defining features, spatial distribution, and motion state. The specific expression of the model is shown in the following equation \mbox{Eq. (\ref{EQ_5})}.

\begin{equation}
\overrightarrow {E_{ji}^v} = \lambda {Q_j}\frac{{{e^{\beta {a_j}\cos \varphi_{j i} }}}}{{\left| {\overrightarrow {{k_{ji}}^\prime } } \right|}^{\beta_2}}\frac{{\overrightarrow {{k_{ji}}^\prime } }}{{\left| {\overrightarrow {{k_{ji}}^\prime } } \right|}},
\label{EQ_5}
\end{equation}
where $\overrightarrow {E_{ji}^v}$ is the interactive vehicle risk field vector on the subject vehicle from an interactive vehicle, $Q_j$ is the virtual energy of the interactive vehicle $j$ defined in \mbox{Eq. (\ref{EQ_2})}, $a_j$ is its acceleration, $\varphi_{j i}$ is the clockwise angle of the subject vehicle $i$ to the velocity direction of the interactive vehicle $j$, ${k_{ji}}^\prime$ is the virtual distance between the interactive vehicle and the subject point $i$, $\lambda$, and $\beta$ and $\beta_2$ are undetermined parameters related to the interactive vehicle risk field strength, the acceleration and its distribution separately.

\begin{figure}[!t]\centering
    \includegraphics[width=14.0cm]{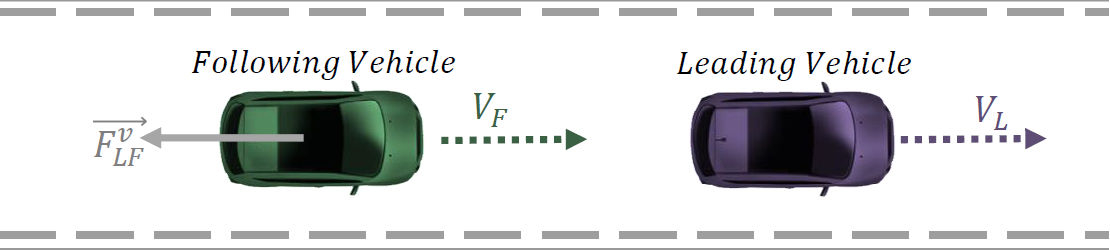}
    \caption{Illustration of the single-lane scenario.}\label{Illustration}
\end{figure}

\mbox{Fig. \ref{Illustration}} shows a typical car-following scenario. Both the leading and following vehicles are traveling along the middle line of the road. The interactive vehicle risk affecting the following vehicle comes from the leading vehicle. In this scenario, the leading vehicle's interactive vehicle risk field strength could be calculated by the equation as \mbox{Eq. (\ref{EQ_6})}.
\begin{equation}
\overrightarrow {E_{LF}^v} = \lambda {Q_L}\frac{{{e^{\beta {a_L}\cos \varphi_{L F} }}}}{{\left| {\overrightarrow {{k_{LF}}^\prime } } \right|}^{\beta_2}}\frac{{\overrightarrow {{k_{LF}}^\prime } }}{{\left| {\overrightarrow {{k_{LF}}^\prime } } \right|}}.
\label{EQ_6}
\end{equation}

The field force of a vehicle in the car-following scenario relates to the degree of risk that the vehicle receives, and the higher the risk, the greater the field force that the vehicle receives. The field force is mainly determined by risk field strength, the vehicle motion state, and vehicle attributes. These factors are summarized as follows:

1. Risk field strength: Field force is a repulsive force and intends to prevent vehicles from approaching other objects, including non-moving obstacles and moving vehicles. Generally, a greater risk field strength indicates a higher risk. Thus the field force received by the vehicle positively correlates with risk field strength.

2. Vehicle Attribute: Attribute is integration of vehicle size and typical motion state. 

1) Vehicle size is a criterion for the classification of vehicle types (e.g., the truck size is usually greater than that of the car). Under the same velocity conditions, the field force of the vehicle with a large size is greater than that of the vehicle with a small size.

2) The vehicle's motion state is typically determined by its velocity. The field force is related to velocity closely. Generally, the vehicle has a bigger field force at a higher velocity than at a lower speed. When the vehicle is stationary, its velocity has no effect on the field force.

Through the above analysis, we can refer to the formula of electric field force: $F = Eq$, in which the field intensity $E$ can be understood as the interactive vehicle risk field strength ${\left| {\overrightarrow {{E_{LF}^v}^\prime } } \right|}$; $q$ is the quantity of electric charge, which belongs to the defined attribute in the electric field environment. In this paper, the defined attribute of the vehicle is its virtual energy $Q_L$ and $Q_F$. Then the force $\overrightarrow {F_{LF}^v}$ of the following vehicle suffered by the leading vehicle can be calculated from \mbox{Eq. (\ref{EQ_7})}.
\begin{equation}
\overrightarrow {F_{LF}^v}  = \overrightarrow {E_L^v} {Q_F} = \lambda {Q_L}{Q_F}\frac{{{e^{\beta {a_L}\cos \varphi_{L F} }}}}{{\left| {\overrightarrow {{k_{LF}}^\prime } } \right|}^{\beta_2}}\frac{{\overrightarrow {{k_{LF}}^\prime } }}{{\left| {\overrightarrow {{k_{LF}}^\prime } } \right|}}.
\label{EQ_7}
\end{equation}

The force $\overrightarrow {F_{LF}^v}$ is a cause of changing the following vehicle's motion state. The formula of Newton's second law, "$F=ma$," is referred to, where the acceleration $a$ is similar to acceleration $\overrightarrow {a_{LF}^v}$ in the proposed model. Through the analysis of the driver's perception above, the mass of an object $m$ can be understood as the defining sizes $L_F^{s_l} W_F^{s_w} $. Therefore, the acceleration $\overrightarrow {a_{LF}^v}$ of the following vehicle under the force $\overrightarrow {F_{LF}^v}$ is shown as \mbox{Eq. (\ref{EQ_8})}.
\begin{equation}
\overrightarrow {a_{LF}^v} = \frac{{\overrightarrow {F_{LF}^v} }}{{L_F^{{s_l}}W_F^{{s_w}}}} = \frac{{\overrightarrow {E_L^v} {Q_F}}}{{L_F^{{s_l}}W_F^{{s_w}}}} = \overrightarrow {E_L^v} {e^{ \alpha {v_F}\cos {\varphi_{F L}}}}
=\lambda L_{L}^{s_{l}} W_{L}^{s_{w}} \frac{e^{\alpha\left(v_{L} \cos {\varphi_{L F}}+v_{F} \cos {\varphi_{F L}}\right)+\beta a_{L} \cos {\varphi_{L F}}}}{\left|\overrightarrow{k_{L F}^{\prime}}\right|^{\beta_{2}}} \frac{\overrightarrow{k_{L F}^{\prime}}}{\left|\overrightarrow{k_{L F}^{\prime}}\right|}.
\label{EQ_8}
\end{equation}

Two basic principles of the acceleration functions encoding the driving behavior in the car-following scenario were proposed in the article \cite{Karaboga:2010}: the acceleration should be a strictly increasing function of the virtual distance; the acceleration is a decreasing function of the speed difference between interactive vehicles.

Firstly,the acceleration is a strictly increasing function of the virtual distance $\overrightarrow{k_{L F}^{\prime}}$ to the leading vehicle. This means the acceleration increases with the virtual distance increases. Thus, \mbox{Eq. (\ref{EQ_k1})} is the following:

\begin{equation}
\frac{\partial(a_{LF}^v(v_{L},v_{F},a_{L},\overrightarrow{k_{L F}^{\prime}}))}{\partial({\left|\overrightarrow{k_{L F}^{\prime}}\right|})}\geq0.
\label{EQ_k1}
\end{equation}

$\overrightarrow{k_{L F}^{\prime}}$ and $\overrightarrow {a_{LF}^v}$ are vectors with the contrary directions. Therefore, $\frac{\overrightarrow{k_{L F}^{\prime}}}{\left|\overrightarrow{k_{L F}^{\prime}}\right|}=-1$.

After simplifying the above formula, the following formula \mbox{Eq. (\ref{EQ_k2})} is obtained:
\begin{equation}
\lambda {\beta_{2}}L_{L}^{s_{l}} W_{L}^{s_{w}}{\left|\overrightarrow{k_{L F}^{\prime}}\right|^{-\beta_{2}-1}} e^{\alpha\left(v_{L} \cos {\varphi_{L F}}+v_{F} \cos {\varphi_{F L}}\right)+\beta a_{L} \cos {\varphi_{L F}}}\geq0.
\label{EQ_k2}
\end{equation}

Meanwhile, $L_{L}$, $W_{L}$, ${\left|\overrightarrow{k_{L F}^{\prime}}\right|}$ and $e^{\alpha\left(v_{L} \cos {\varphi_{L F}}+v_{F} \cos {\varphi_{F L}}\right)+\beta a_{L} \cos {\varphi_{L F}}}$ are greater than 0. Therefore, $\lambda {\beta_{2}}$ is required to be greater than 0 in calibration.

Secondly, the acceleration is a decreasing function of the speed difference between the subject vehicle and the interactive vehicle in the car-following scenario. For example, deceleration needs to be enhanced when the following vehicle's velocity is much faster than that of the leading vehicle. And acceleration increases when the speed of the following vehicle is much slower than that of the leading vehicle. We order $dv=v_{F}-v_{L}$. Thus,
\begin{equation}
\overrightarrow {a_{LF}^v}=\lambda L_{L}^{s_{l}} W_{L}^{s_{w}} \frac{e^{\alpha{v_{L}} \left(\cos {\varphi_{L F}}+ \cos {\varphi_{F L}}\right)+\alpha dv\cos {\varphi_{F L}}+\beta a_{L} \cos {\varphi_{L F}}}}{\left|\overrightarrow{k_{L F}^{\prime}}\right|^{\beta_{2}}} \frac{\overrightarrow{k_{L F}^{\prime}}}{\left|\overrightarrow{k_{L F}^{\prime}}\right|},
\label{EQ_dv1}
\end{equation}

\begin{equation}
\frac{\partial(a_{LF}^v(v_{L},dv,a_{L},\overrightarrow{k_{L F}^{\prime}}))}{\partial(dv)}\leq0.
\label{EQ_dv2}
\end{equation}

Then, the following formula \mbox{Eq. (\ref{EQ_dv3})} is required:
\begin{equation}
-\lambda \alpha \cos {\varphi_{F L}} L_{L}^{s_{l}} W_{L}^{s_{w}} \left|\overrightarrow{k_{LF}^{\prime}}\right|^{-\beta_{2}} {e^{\alpha{v_{L}} \left(\cos {\varphi_{L F}}+ \cos {\varphi_{F L}}\right)+\alpha dv\cos {\varphi_{F L}}+\beta a_{L} \cos {\varphi_{L F}}}} \leq0.
\label{EQ_dv3}
\end{equation}

Meanwhile, $L_{L}$, $W_{L}$, ${\left|\overrightarrow{k_{L F}^{\prime}}\right|}$ and ${e^{\alpha{v_{L}} \left(\cos {\varphi_{L F}}+ \cos {\varphi_{F L}}\right)+\alpha dv\cos {\varphi_{F L}}+\beta a_{L} \cos {\varphi_{L F}}}}$ are greater than 0. $\cos {\varphi_{F L}}$ is also positive because vehicle $L$ is in front of vehicle $F$. 
Therefore, $\lambda {\alpha}$ is also required to be greater than 0 in calibration.

\subsection{Restrictions Risk Surrogate}

As mentioned in subsection 2.1, the restrictions risk surrogate characterizes the physical field of the impact of impassable areas on driving safety. The attributes of restriction mainly determine the field strength and vector direction of the restrictions risk surrogate.

Restrictions forming the risk surrogate can be divided into two categories.

1. Stationary objects that could substantially collide with vehicles and cause losses (personal and property), such as isolation, anti-collision pier, and construction sign. The closer the vehicle is to this type of stationary object, the more likely it is to collide with it. Moreover, as the vehicle approaches, this possibility and the risk of collisions do not increase linearly since as the clearance decreases, the driving risk increases faster. Meanwhile, as the vehicle approaches this kind of stationary object, the variation law of the degree of driving risk posed by this kind of stationary object is the same no matter from any direction. 

2. Static objects constrain the behavior of the subject vehicle, such as lane lines, road signs, etc. For this category of restrictions, the violation of the restrictions generated by the vehicle will not lead to a substantial collision accident. However, it will result in a violation of laws or traffic regulations. The driving risk caused by such violations is determined by the traffic regulations they violate. For example, China's traffic regulations stipulate that: vehicles are prohibited from crossing the solid white line, changing lanes, or driving on the line. If a driver finds himself deviating from his current lane, the driver perceives a risk of violating the lane line. This risk drives the vehicle back toward the centerline. In addition, the closer the vehicle is to the lane line, the greater the driving risk.

Then, the undetermined coefficients ${T_\omega }$ that characterize the field strength can be defined according to the degree of different restrictions. For example, when the vehicle is changing lanes, ${T_\omega }$ of the white dotted line is equal to 0, so that the white dotted line has no restriction on a lane-changing vehicle. The restrictions risk field strength ${\left| {\overrightarrow {{E_{i}^t}^\prime } } \right|}$ could be expressed specifically as follows \mbox{Eq. (\ref{EQ_9})}.
\begin{equation}
\overrightarrow {E_{\omega i}^t} = {{T_\omega }\left( {\frac{1}{{\left| {\overrightarrow {{k_{\omega i}}^\prime } } \right| - \left| {\overrightarrow {{k_{dim}}^\prime } } \right|}}} \right)} ^2\frac{{\overrightarrow {{k_{\omega i }}} }}{{\left| {\overrightarrow {{k_{\omega i }}} } \right|}},
\label{EQ_9}
\end{equation}
where, $\overrightarrow {{k_{\omega i }}}$ represents the vector distance from the vehicle $i$ to the restriction $\omega$, $\overrightarrow {{k_{dim}}}$ represents the vector distance between the mass center of the vehicle $i$ to the edge of the vehicle body in the same direction as $\overrightarrow {{k_{\omega i}}}$.

The formula of electric field force: $F=Eq$ is also used, in which the field intensity $E$ can be understood as the restrictions risk field vector $\overrightarrow {E_{\omega i}^t}$; Then the restriction force $\overrightarrow {F_{\omega i}^t}$ formed by restriction risks can be calculated in \mbox{Eq. (\ref{EQ_10})}.
\begin{equation}
\overrightarrow {F_{\omega i}^t} ={Q_i} \overrightarrow {E_{\omega i}^t} = {Q_i}{{T_\omega }\left( {\frac{1}{{\left| {\overrightarrow {{k_{{\omega i} }}^\prime } } \right| - \left| {\overrightarrow {{k_{dim}}^\prime } } \right|}}} \right)} ^2\frac{{\overrightarrow {{k_{{\omega i} }}} }}{{\left| {\overrightarrow {{k_{\omega i}}} } \right|}}.
\label{EQ_10}
\end{equation}

The restriction forces $\overrightarrow {F_{\omega i}^t}$ is the cause of the related motion state to the restriction object $\omega$. The formula of Newton's second law: $F=ma$ is also referred to, where the acceleration $a$ is similar to acceleration $\overrightarrow {a_{\omega i}^t}$. Therefore, the acceleration $\overrightarrow {a_{\omega i}^t}$ of the following vehicle under the action of the restriction force $\overrightarrow {F_{\omega i}^t}$ is shown as \mbox{Eq. (\ref{EQ_11})}.
\begin{equation}
\overrightarrow {a_{\omega i}^t} = \frac{{\overrightarrow {F_{\omega i}^t} }}{{L_i^{{s_l}}W_i^{{s_w}}}} = \frac{{\overrightarrow {E_{\omega i}^t} {Q_i}}}{{L_i^{{s_l}}W_i^{{s_w}}}} = \overrightarrow {E_{\omega i}^t} {e^{- \alpha {v_i}\cos {\varphi _{\omega i}}}}.
\label{EQ_11}
\end{equation}

The acceleration is a strictly increasing function of the distance $\overrightarrow {E_{\omega i}}$ to the restriction object. This means that the acceleration increases with the distance increases. Thus, \mbox{Eq. (\ref{EQ_w1})} is following:

\begin{equation}
\frac{\partial(a_{\omega i}^t(v_{i},\overrightarrow {E_{\omega i}}))}{\partial({\left|\overrightarrow {E_{\omega i}}\right|})}\geq0.
\label{EQ_w1}
\end{equation}

${\overrightarrow {{k_{{\omega i} }}} }$ and $\overrightarrow {a_{\omega i}^t}$ are vectors with the contrary directions. Generally, $\frac{{\overrightarrow {{k_{{\omega i} }}} }}{{\left| {\overrightarrow {{k_{\omega i}}} } \right|}}=-1$.

After simplifying \mbox{Eq. (\ref{EQ_w1})}, the following formula is obtained:
\begin{equation}
 2{{T_\omega }\left( {\frac{1}{{\left| {\overrightarrow {{k_{\omega i}}^\prime } } \right| - \left| {\overrightarrow {{k_{dim}}^\prime } } \right|}}} \right)} ^3{e^{- \alpha {v_i}\cos {\varphi _{\omega i}}}}\geq0.
\label{EQ_w2}
\end{equation}

Meanwhile, $\left| {\overrightarrow {{k_{\omega i}}^\prime } } \right|$, $\left| {\overrightarrow {{k_{dim}}^\prime } } \right|$ and ${e^{- \alpha {v_i}\cos {\varphi _{\omega i}}}}$ are greater than 0. Therefore, undetermined parameters ${T_\omega }$ is required to be greater than 0 in calibration.

\subsection{Speed Risk Surrogate}

According to the previous description of the driving risk, the field force formed by DRS belongs to the short-range force. This means the force will only work on the vehicle under a certain distance. When the vehicles are far apart, their motion behavior is mainly affected by their desired velocity. If the desired velocity is not reached, the vehicle will accelerate. Otherwise, the vehicle will decelerate when its velocity exceeds the desired velocity. Therefore, a velocity potential field can be defined to describe such driving behavior accurately.

The desired velocity of the object vehicle would be affected by the undetermined flow speed $v_{0}$ and speed $v_{l}$ of its interaction vehicle.
\begin{equation}
v_{i}^{d}=\gamma v_{0}+(1-\gamma) v_{l}, \\
\label{EQ_12}
\end{equation}
where, $v_i^d$ is the desired velocity and $\gamma$ is undetermined parameters related to the weight.

The field function named “speed risk surrogate” is mainly affected by the relationship between the desired velocity and the current velocity. It has the following expression, as shown in \mbox{Eq. (\ref{EQ_13})}.
\begin{equation}
\begin{gathered}
\overrightarrow{E_{i}^{s}}=L_{i}^{s_{l}} W_{i}^{s_{w}} {E_{\max }}\left(\frac{|\Delta v|}{v_{i}^{d}}\right)^{\sigma^\prime}\frac{{\overrightarrow {{v_{i}}} }}{{\left| {\overrightarrow {{v_{i}}} } \right|}}, \\
\label{EQ_13}
\end{gathered}
\end{equation}
where, $\overrightarrow {E_i^s}$ is the speed risk field vector, ${E_{\max}}$ represents the maximum field strength, $\Delta v$ is defined by $(v_i^d-v_i)$, $v_i^d$ is the desired velocity, $v_i$ represents the velocity of vehicle $i$, and $\sigma^\prime$ is an undetermined coefficient.

In the speed risk surrogate, which is different from the interactive vehicle risk surrogate and the restrictions risk surrogate, the speed risk force is related to the gradient of the speed risk surrogate. The following equation shows the speed risk force $\overrightarrow {F_{i}^s}$. 
\begin{equation}
\begin{gathered}
\overrightarrow{F_{i}^{s}}=\frac{\mathrm{d}\overrightarrow{{E_{i}^{s}}}}{\mathrm{d}{\Delta v}} \frac{{\overrightarrow {{v_{i}}} }}{{\left| {\overrightarrow {{v_{i}}} } \right|}}.\\
\label{EQ_14}
\end{gathered}
\end{equation}

Meanwhile, according to the formula of Newton's second law, the acceleration $\overrightarrow {a_i^s}$ of the subject vehicle under the speed risk force is as follows \mbox{Eq. (\ref{EQ_15})}.
\begin{equation}
\begin{gathered}
\overrightarrow{F_{i}^{s}}=L_{i}^{s_{l}} W_{i}^{s_{w}} \overrightarrow{a_{i}^{s}}. \\
\label{EQ_15}
\end{gathered}
\end{equation}

Thus the acceleration is mainly described by the relationship between the current velocity and desired velocity of the subject vehicle. The following equation \mbox{Eq. (\ref{EQ_16})} is derived from \mbox{Eq. (\ref{EQ_14})} and \mbox{Eq. (\ref{EQ_15})}.
\begin{equation}
\begin{gathered}
\overrightarrow{a_{i}^{s}}= \frac{\overrightarrow{F_{i}^{s}}}{L_{i}^{s_{l}} W_{i}^{s_{w}}} = \frac{1}{L_{i}^{s_{l}} W_{i}^{s_{w}}}  \frac{\mathrm{d}\overrightarrow{{E_{i}^{s}}}}{\mathrm{d}{\Delta v}} = \left\{\begin{array}{l}
{E_{\max }} \frac{\sigma^{\prime}}{{v_{i}^{d}}^{\sigma^{\prime}}} {|\Delta v|}^{\sigma^{\prime}-1} \frac{{\overrightarrow {{v_{i}}} }}{{\left| {\overrightarrow {{v_{i}}} } \right|}} , \Delta v \geq 0. \\
-{E_{\max }} \frac{\sigma^{\prime}}{{v_{i}^{d}}^{\sigma^{\prime}}} {|\Delta v|}^{\sigma^{\prime}-1} \frac{{\overrightarrow {{v_{i}}} }}{{\left| {\overrightarrow {{v_{i}}} } \right|}}, \Delta v<0.
\end{array}\right.\\
\label{EQ_16}
\end{gathered}
\end{equation}

Furthermore, we order ${a_{\max }}={E_{\max }} \frac{\sigma^{\prime}}{{v_{i}^{d}}^{\sigma^{\prime}}}$ and $\sigma=\sigma^{\prime}-1$.Then, the simplified formula is \mbox{Eq. (\ref{EQ_17})}.
\begin{equation}
\begin{gathered}
\overrightarrow{a_{i}^{s}}=\left\{\begin{array}{l}
{a_{\max }}{|\Delta v|}^{\sigma}\frac{{\overrightarrow {{v_{i}}} }}{{\left| {\overrightarrow {{v_{i}}} } \right|}} , \Delta v \geq 0. \\
-{a_{\max }}{|\Delta v|}^{\sigma}\frac{{\overrightarrow {{v_{i}}} }}{{\left| {\overrightarrow {{v_{i}}} } \right|}} , \Delta v<0.
\end{array}\right. \\
\label{EQ_17}
\end{gathered}
\end{equation}
where $\overrightarrow {a_i^s}$ represents the acceleration vector of a vehicle $i$ affected by its desired speed. The formula depicts that when $\Delta v \ge 0$, the acceleration is positive, in the next moment, the vehicle is in an accelerated state; when $\Delta v < 0$, the acceleration is negative, in the next moment, the vehicle is in a decelerated state. 

The acceleration is a decreasing function of the desired velocity ${v_{i}^{d}}$. Moreover, the vehicle accelerates towards the desired velocity if not constrained by other vehicles or restrictions:

\begin{equation}
\frac{\mathrm{d}a_{i}^{s}}{\mathrm{d}v_{i}}<0.
\label{EQ_vd1}
\end{equation}

$\overrightarrow {{v_{i}}}$ and $\overrightarrow{a_{i}^{s}}$ are motion vectors with the same initial direction. Generally, $\frac{{\overrightarrow {{v_{i}}} }}{{\left| {\overrightarrow {{v_{i}}} } \right|}}=1$.

When $\Delta v={v_i^d-v_i}\geq 0$, the following formula is derived:
\begin{equation}
\frac{\mathrm{d}a_{i}^{s}}{\mathrm{d}v_{i}}=
-{E_{\max }}\frac{\sigma^{\prime}({\sigma^{\prime}-1})} {{v_{i}^{d}}^{\sigma^{\prime}}} \left(v_i^d-v_i\right)^{\sigma^{\prime}-2}<0.
\label{EQ_vd2}
\end{equation}

When $\Delta v={v_i^d-v_i}<0$, the following formula is derived:
\begin{equation}
\frac{\mathrm{d}a_{i}^{s}}{\mathrm{d}v_{i}}=
-{E_{\max }}\frac{\sigma^{\prime}({\sigma^{\prime}-1})} {{v_{i}^{d}}^{\sigma^{\prime}}} \left(v_i-v_i^d\right)^{\sigma^{\prime}-2}<0.
\label{EQ_vd3}
\end{equation}

Furthermore, ${E_{\max }}$ and $v_i^d$ are greater than 0, as the same as ${|\Delta v|}^{\sigma^{\prime}-2}$. Therefore, ${\sigma^{\prime}({\sigma^{\prime}-1})}$ is required to be greater than 0. Due to $\sigma=\sigma^{\prime}-1$, $\sigma>0$ or $\sigma<-1$ are limits of $\sigma$ in calibration.

\begin{figure}
\centering
\includegraphics[width=1.0\textwidth]{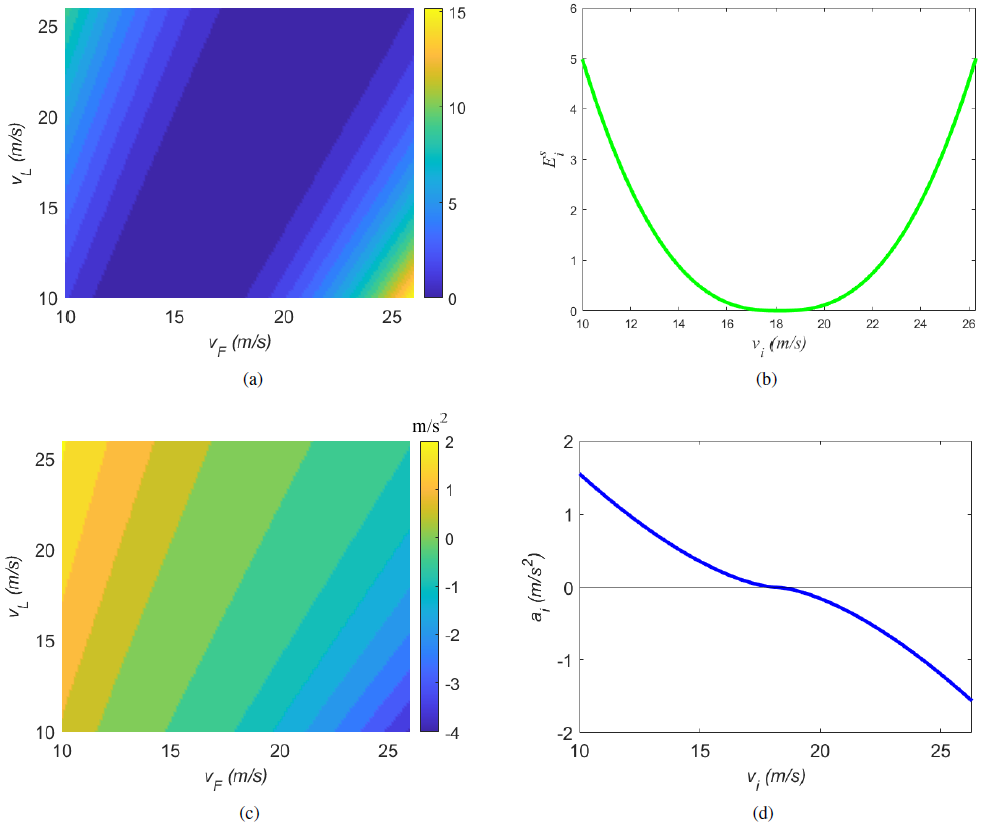}
\caption{\label{Visualization of risk reflection}Distribution map of risk field strength (yellow indicates higher risk, while blue is the opposite) and acceleration (yellow indicates acceleration, while blue is deceleration) from speed risk surrogate.}
\end{figure}

Based on the establishment of speed risk surrogate, figures about risk field strength and acceleration are shown in \mbox{Fig. \ref{speed risk surrogate}}. \mbox{Fig. \ref{speed risk surrogate}}(a) displays the map of the risk field strength (yellow indicates higher risk, while blue is the opposite) in different velocity scenarios of car-following. The X-axis is the velocity of the following vehicle, and the Y-axis is that of its leading vehicle. The speed risk field strength increases monotonically and nonlinearly with the absolute speed difference. The speed risk energy is not lowest when the two vehicles are at the same velocity because of the influence of the speed of traffic flow. Specifically, \mbox{Fig. \ref{speed risk surrogate}}(b) shows the risk curve of different velocity plans when the interactive velocity is driving at 17m/s. The lowest speed risk field strength does not occur at 17m/s. \mbox{Fig. \ref{speed risk surrogate}}(c) displays the map of the acceleration (yellow indicates acceleration, while blue is deceleration) in different velocity scenarios of car-following. As the same as \mbox{Fig. \ref{speed risk surrogate}}(a), the X-axis is the velocity of the following vehicle, and the Y-axis is that of its leading vehicle. The following vehicle will accelerate when its velocity is lower than that of its desired velocity while decelerating when its velocity is higher. \mbox{Fig. \ref{speed risk surrogate}}(b) shows the acceleration curve of different velocity plans when the interactive velocity is driving at 17m/s. This clearly shows the relationship between acceleration and desired velocity. The gradient of the acceleration curve decreases with the absolute velocity difference. If the vehicle speed is much slower than the desired speed, the acceleration will increase obviously, which is in line with the reality of drivers' pursuit of driving efficiency. On the contrary, when the vehicle speed is much faster than the ideal vehicle speed, the deceleration increases to achieve driving safety.

\section{Calibration}
Model calibration is a necessary process to concretize the undetermined parameters. The parameter calibration process is completed based on meta-heuristic algorithms and a naturalistic driving dataset. We define many undetermined parameters in section 2 and those parameters are of great value for the proposed model. In this section, we intend to select driving data and concretize those parameters effectively to improve the model's applicability and prevent its distortion.

\subsection{Data Selection and Pre-processing}

Safe interactions generally happen in naturalistic driving behavior. Although human drivers cannot obtain the specific movement information of surrounding vehicles, the human driver has the most sophisticated task processor in the world and can sensitively perceive the driving risks of the surrounding environment. Especially when the speed of the front vehicle changes obviously, the following manual-driven vehicle as the rear vehicle will change the speed accordingly. Therefore, the real driving trajectory dataset can be used as the research dataset of this paper and be adopted to calibrate the proposed model. In order to model DRS more in line with Chinese human driving habits, this paper chooses a dataset of naturalistic driving from Ubiquitous Traffic Eyes\cite{LiZhibin}  in Nanjing City, China, to study the model performance.

\begin{figure}[!t]\centering
    \includegraphics[width=14.0cm]{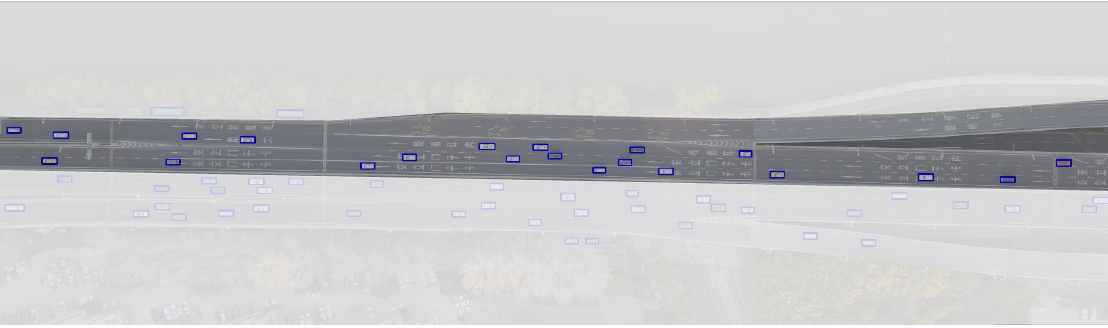}
    \caption{The road segment of the naturalistic driving dataset.}\label{fig_4}
\end{figure}

The dataset contains pure trajectory data, which applies to micro-driving behavior studies. According to the purpose of this study, we select the trajectory data on a segment of an urban highway that is one of the busiest interchange hubs in the south of Nanjing city, China, as shown in \mbox{Fig. \ref{fig_4}}. A large number of naturalistic driving activities occur in this area. The length of the road section is 352 meters, east-west direction. The dataset is captured by unmanned aerial vehicles at 285m and provid the precise vehicle position coordinates with a time accuracy of 0.1s. The detailed data of each vehicle in 1184 seconds are contained, including vehicle ID, position, speed, acceleration, width, length, and time stamp.

Before using the dataset for our study, it is necessary to pre-process the original data. First, all data at the same time are taken out. They are grouped according to lane lines and then sort by X-axis coordinates within each lane to realize the judgment of the leading and following vehicles. The track pairs are searched according to vehicle IDs. Finally, we filter out the track segments where the lane-changing behavior occurs less than 0.5s.
Based on the data pre-processing, 3063 pairs of car-following trajectories are obtained for further verification.

\subsection{Model Calibration}
Parameter calibration is a bridge from the theoretical level to the practical application level. An effective parameter calibration process can enhance the model's applicability and avoid distortion. In this paper, a modified particle swarm optimization algorithm is selected to concretize the undetermined parameters in the model combined with naturalistic driving data. The Root Mean Squared Error (RMSE) is adopted as the loss function, which describes the deviation between the real value and the simulation value and is defined as an index for evaluating and verifying parameter calibration results. The specific expression of RMSE is shown in the following equation \mbox{Eq. (\ref{EQ_18})}.
\begin{equation}
\begin{gathered}
RMSE=\sqrt{\frac{\sum_{i=1}^{N}\left(x_{i}^{\text {real }}-x_{i}^{\text {sim }}\right)^{2}}{N}}.
\label{EQ_18}
\end{gathered}
\end{equation}
where, $x_i^\text{real}$ represents the actual position of sample $i$, $x_i^\text{sim}$ represents its simulation position, and $N$ is the sample size.

In the design of the particle swarm optimization algorithm, to avoid falling into local optimization while iteratively calculating the local optimal particles in the particle pool, the standby particle pool is introduced for iterative calculation, and the better particles are replaced by the particles in the previous particle pool. Based on the calibration, the proposed model parameters are obtained and shown in \mbox{Table \ref{table_1}}. Meanwhile, two models, IDM and DRPFM(a risk field-based model for car-following) \cite{li2020dynamic}, are selected as the baseline for comparison purposes of the proposed model, and the same optimization algorithm is applied to calibrate their parameters which are shown in \mbox{Table \ref{table_1}}.

\begin{table}[]
\caption{Parameters calibration of IDM, DRPFM, amd DRS.}\label{table_1}
\resizebox{\textwidth}{!}{%
\begin{tabular}{lllllllllllll}
\hline
\multicolumn{13}{c}{IDM}                                                                           \\ \hline
Parameter & $v_{free}$          & $a_{max}$     & $b$      & $s_{0}$     & $T$      & $\delta$ &   &   &   &    &    &    \\
Value & 21.9612 & 0.4708 & 1.1153  & 3.0000 & 1.5754  & 4.0000 &        &       &       &       &       &       \\ \hline
\multicolumn{13}{c}{DRPFM}                                                                         \\ \hline
Parameter & $v_{0}^{(a)}$     & $a_{max}$     & $\alpha$      & $\lambda$      & $\sigma$    & $\beta_{1}$      &$\beta_{2}$ & &   &    &    &    \\
Value & 27.7611 & 0.5083 & -0.3627 & 0.0096 & 0.2538  & -0.6714 & 0.4502 &  &       &       &       &       \\ \hline
\multicolumn{13}{c}{DRS}                                                       \\ \hline
Parameter & $v_{0}$          & $\alpha$     & $\beta$      & $\lambda$      & $\gamma$     & $\sigma$     & $a_{\max }$      & $\beta_{2}$     & $s_l$    &$s_w$     &$w_{1}$     & $w_{2}$  \\
Value & 18.220 & 5.425 & 0.185  & 0.374 & 14.000 & 0.576 & 1.554  & 0.984 & 0.068 & 0.007 & 0.239 & 0.881 \\ \hline
\end{tabular}%
}
\end{table}

\section{Simulation and Discussion}
In the analysis of safe interactions, more safety metrics do not work because safety is the precondition. Meanwhile, the safety metrics are defined subjectively, and there is no unified threshold to claim. This is why we derived the acceleration from field strength in above \mbox{Eq. (\ref{EQ_8})}, \mbox{Eq. (\ref{EQ_11})}, and \mbox{Eq. (\ref{EQ_17})}. The acceleration makes the defined safety metric relate to an object parameter. Subsection 4.1 shows the acceleration variations which indicate safety in different car-following scenarios. Trajectory pairs of vehicles are typical features of human driving interactions. In subsection 4.2, trajectory estimation derived from acceleration is presented to show the results of safe interactions. Compared with other models, the proposed model show a better performance in safe interactions using the real-world dataset. To further discuss detailed performance and models' characteristics, key factors of the safe interactions processes are shown in subsection 4.3, such as acceleration, space headway, space headway error, and velocity. The discussions focus on safe interaction features are presented according to their performances.

\subsection{Visualization of Risk Reflection Variation}
In order to better reflect the ability of risk assessment, simulation experiments are conducted for IDM, DRPFM, and DRS. As we all know, TTC is time-based and is widely used to characterize risk. IDM and DRPFM are acceleration-based, and their accelerations can be used as reflections of the driving risk. DRS is a risk assessment model, and there is heterogeneity in calibrating the threshold value of risk degree because the threshold value may change with different motion scenarios. Since we derived acceleration from the field strength in the proposed model. The acceleration is selected as the response to risk degree. This provides a parameter for model comparison and solves the problem of inconsistent thresholds.

In the simulation, the leading vehicle is driving forward at a constant speed of 20m/s. The following vehicle is driving behind over different speeds, and different gaps and the related risk reflection variations of the four algorithms are shown in \mbox{Fig. \ref{Visualization of risk reflection}}. The gap represented the distance between the leading vehicle and the following vehicle. It is selected as the risk reflection of TTC, while acceleration is used for IDM, DRPFM, and DRS. 

When the following vehicle is slower than the leading vehicle, TTC shows absolute safety. TTC shows the remaining operation time and reflects the risk variation when the following vehicle is higher than the leading vehicle. It can be seen that the risk response of TTC is consistent with the actual situation. As we all know, deceleration is impossible to get $10m/s^2$ due to the limitation of vehicle kinetics. As the velocity of the following vehicle increases, IDM and DRPFM achieve the unattainable deceleration $10m/s^2$ earlier than DRS, even when the situation is safe in TTC. Meanwhile, the following vehicle does not accelerate when the velocity is much slower than that of the leading vehicle in IDM and DRPFM, while DRS does. Further, IDM and DRPFM have a more fierce transition from "blue" to "green" than DRS. For vehicles in motion, the speed and space headway vary continuously. The homogeneity of DRS risk response is better for the distribution of risk response under different speeds and gaps. That is, its jerk of acceleration change is smaller. However, the acceleration change gradient of IDM and DRPFM is larger, resulting in greater jerk. When the vehicle accelerates, passengers will feel uncomfortable. This discomfort comes not only from acceleration but also from a jerk. Therefore, the DRS with a smaller jerk can bring a more comfortable ride. It is worth noting that this section shows the risk reflection variation of the proposed model compared with TTC and two car-following models. Numerical simulations show risk variation and distribution trends rather than random realizations. If surrounding vehicles change motion status, the specific values will vary while the trend and distribution will not.

\begin{figure}
\centering
\includegraphics[width=1.0\textwidth]{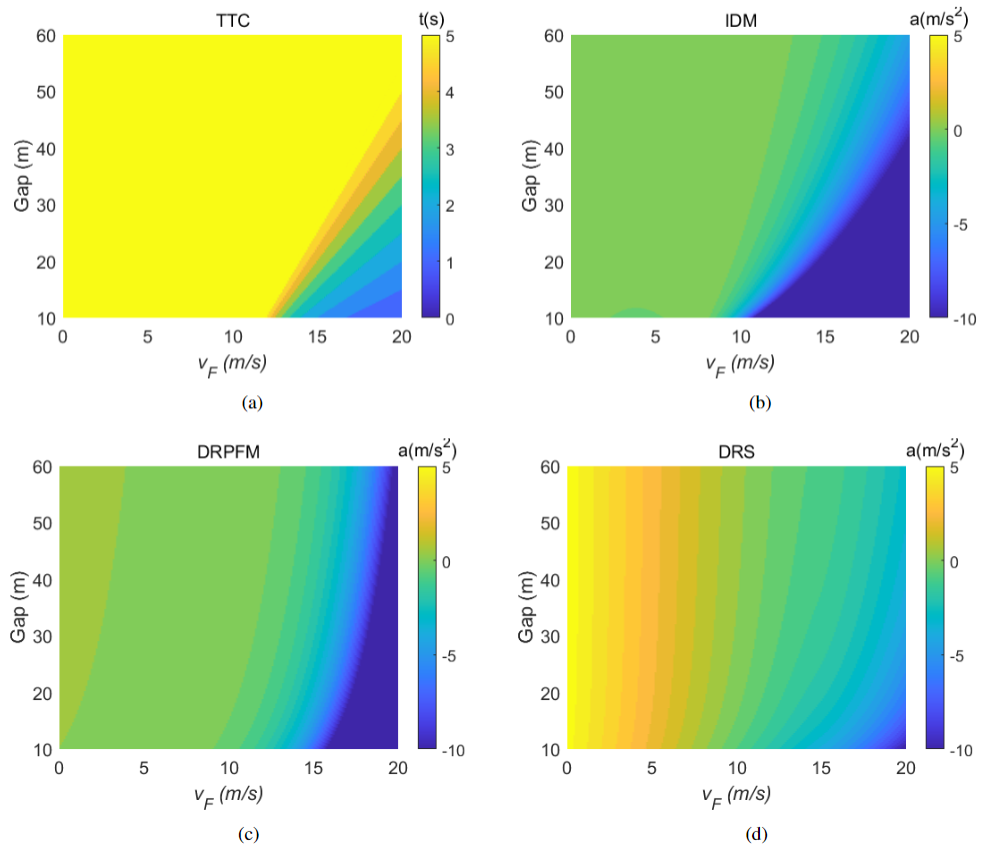}
\caption{\label{Visualization of risk reflection}Visualization of risk reflection variation using different models (blue indicates higher risk, while yellow is the opposite).}
\end{figure}

\subsection{Trajectory Estimation Performance}
\begin{figure}
\centering
\includegraphics[width=1.0\textwidth]{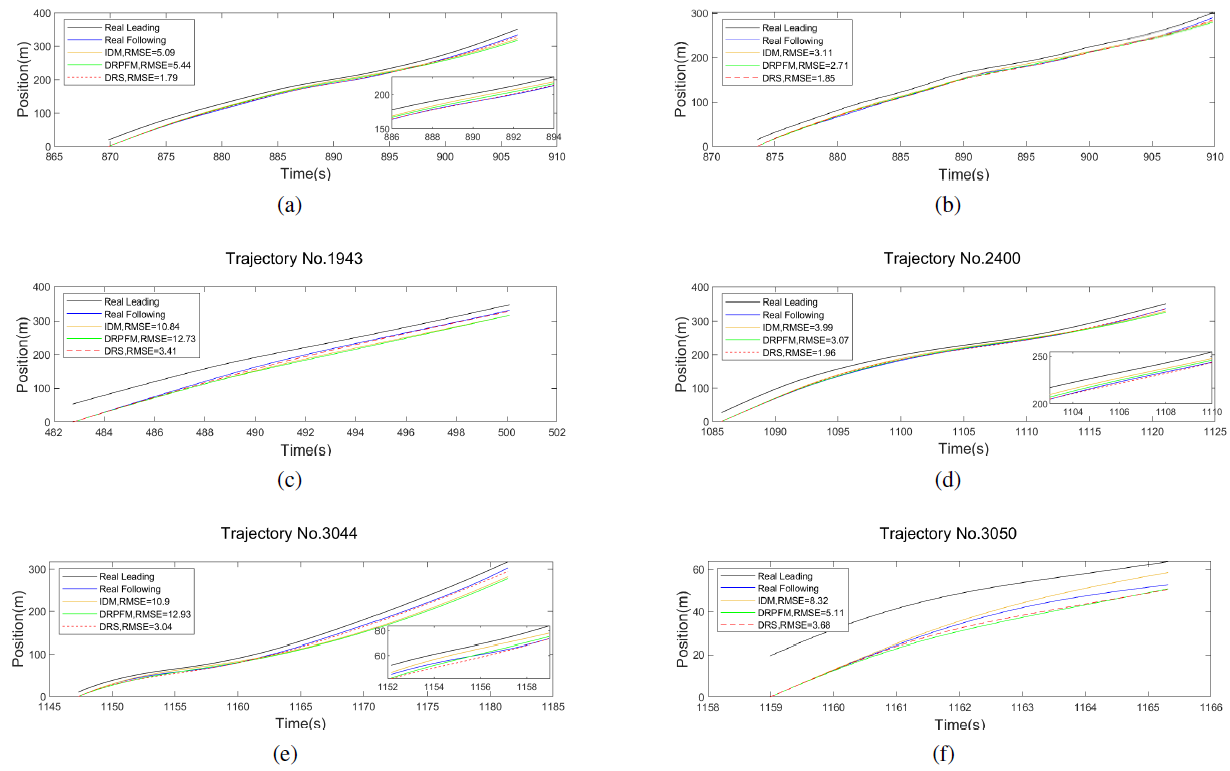}
\caption{\label{Trajectory estimation}Trajectory estimation performance of IDM, DRPFM and DRS.}
\end{figure}

Trajectory pairs of vehicles are typical features of human driving interactions. To discuss the performances on safe interaction, simulations of trajectory estimation are completed based on a real-world dataset. The same naturalistic driving dataset and optimization algorithm are used to calibrate IDM, DRPFM, and DRS parameters, respectively. The simulation results of the rear vehicle are compared with the real trajectory, and the results are shown in \mbox{Fig. \ref{Trajectory estimation}}. The black line in the figure shows the real track of the front vehicle, the blue line is the real track of the rear vehicle, and the yellow and green are the trajectory estimation of the rear vehicle obtained by IDM and DRPFM. The red dotted line is the estimation result of the rear vehicle trajectory under the proposed method DRS. Legend also shows the RMSEs of the three algorithms.

The comparison shows that the three models all perform well with no collision, which shows as lines crossing. Phenomena appear that in IDM, the vehicle is driving too close to the front vehicle in \mbox{Fig. \ref{Trajectory estimation}} (a), (d) and (e). We use RMSE to evaluate the trajectory estimation performance, shown in legends. Although IDM, DRPFM, and DRS cannot completely match the real data, their RMSE of trajectory estimation is controlled within a reasonable range in general. DRS usually maintains a smaller RMSE than IDM and DRPFM in driving interactions in different car-following cases. 

\begin{figure}[!t]\centering
    \includegraphics[width=12.0cm]{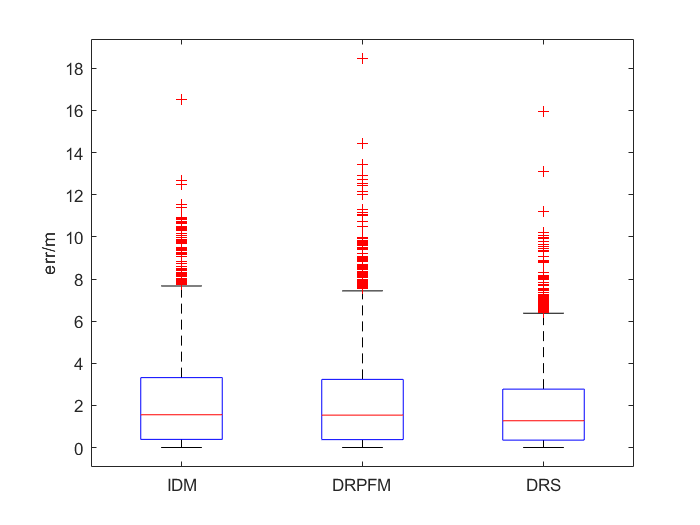}
    \caption{Boxplot for trajectory estimation RMSE of IDM, DRPFM and DRS.}\label{Boxplot}
\end{figure}

To avoid the particularity caused by individual trajectory analysis, average RMSEs of statistics are calculated and show that our model, DRS (3.334), performs most accurately compared with IDM (3.9507) and DRPFM (3.9392). Further, to compare the performance in trajectory estimation of different models more convincingly and comprehensively, RMSEs obtained by simulating all trajectories are drawn into the boxplot, as shown in \mbox{Fig. \ref{Boxplot}}. The boxplot can reflect the dependability of the trajectory estimation performance, and the interquartile range (IQR) is used to present a discrete degree of value.

As can be seen in \mbox{Fig. \ref{Boxplot}}, the median of the boxplot produced by the DRS algorithm is 1.2698, and other medians of the boxplot produced are 1.5563 (IDM) and 1.5393 (DRPFM). The three models have the same lower adjacent, and DRS has the lowest upper adjacent, which indicates that the overall RMSE of trajectory estimation by DRS is lower than\ other algorithms. Besides, the IOR of the boxplot generated by the DRS algorithm is 2.42, and the IOR of the boxplot generated by the other algorithms is 2.93 (IDM) and 2.86 (DRPFM), respectively. This indicates that DRS's RMSE varies less than IDM and DRPFM. The more consistent RMSE should make estimations more dependable than others. Based on the above 
characteristic, DRS is more dependable, and its RMSEs remain at lower levels, while IDM and DRPFM are more variable and at higher RMSE levels. This proves that the DRS algorithm performs better in car-following interactions than the IDM and DRPFM algorithms. 

\subsection{Detail Performance in Velocity, Headway and Acceleration}
\begin{figure}
\centering
\includegraphics[width=1.0\textwidth]{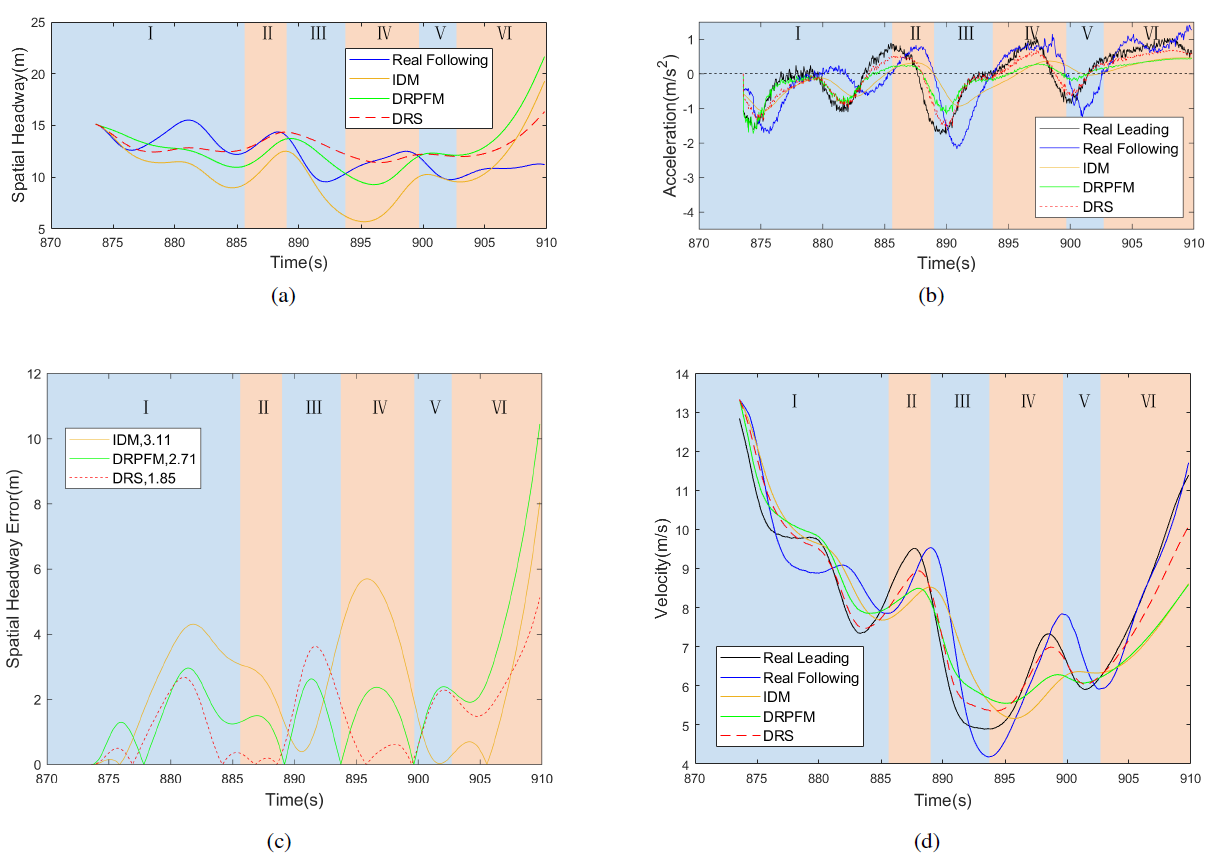}
\caption{\label{1vha}The comparison of performance about velocity, spatial headway and acceleration in trajectory No.1396.}
\end{figure}

\begin{figure}
\centering
\includegraphics[width=1.0\textwidth]{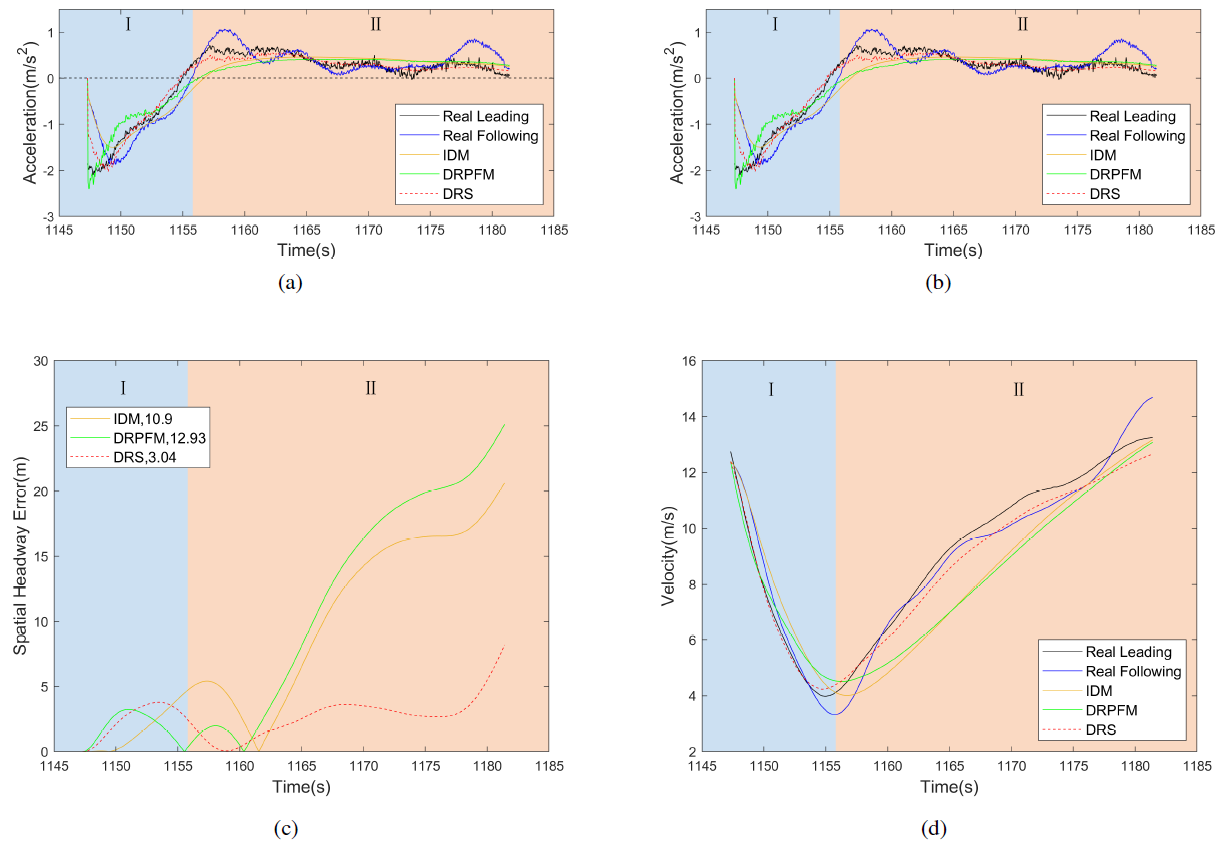}
\caption{\label{2vha}The comparison of performance about velocity, spatial headway and acceleration in trajectory No.3044.}
\end{figure}

To further discuss detailed performance in safe interactions, acceleration, space headway, space headway error, and velocity are selected as key factors. This section chose two cases from trajectory pairs in the naturalistic driving dataset and discussed the performance in the safe interaction process. \mbox{Fig. \ref{1vha}} and \mbox{Fig. \ref{2vha}} show the simulation results of the three algorithms in velocity, space headway, error, and acceleration. The black and blue solid lines in the figure represent the real data of the front and rear vehicles. Yellow, green, and red lines are IDM, DRPFM, and DRS simulation results.
To facilitate subsequent analysis, acceleration scenarios and deceleration scenarios of the following vehicle are drawn in the orange block and blue block, respectively, in each figure, and Roman numerals are marked on each area.

\mbox{Fig. \ref{1vha}} (a) and \mbox{Fig. \ref{2vha}} (a) show space headway of trajectory No.1396 and No.3044. When the leading vehicle decelerates ( \uppercase\expandafter{\romannumeral1}, \uppercase\expandafter{\romannumeral3}, \uppercase\expandafter{\romannumeral5}, \uppercase\expandafter{\romannumeral7}), the space headway curve of IDM is the lowest, while DRPFM's and DRS's are higher, indicating that the IDM model is radical and DRPFM and DRS are more conservative. Space headways in DRPFM and DRS are not too close to the front vehicle due to their risk assessment considering velocity and acceleration, while IDM does not consider acceleration. As shown in area \uppercase\expandafter{\romannumeral3} and \uppercase\expandafter{\romannumeral5} in \mbox{Fig. \ref{1vha}}, with the sharp deceleration of the vehicle in front, IDM takes into account the reduction of the speed of the leading vehicle and sharply reduces the following gap. However, in this case, the deceleration of the leading vehicle is very large, even reaching $2m/s^2$, and a risk of collision may occur. In the same scenario, DRPFM and DRS consider the impact of the vehicle's deceleration in front and reserv a longer following distance for driving safety.

Meanwhile, the proposed DRS can better fit the real space headway when vehicles speed up. Further, \mbox{Fig. \ref{1vha}} (b) and \mbox{Fig. \ref{2vha}} (b) show RMSE curves of headway estimations. The mean RMSEs of the three models are also shown on the top left corner of each graph, and the DRS has a minimum mean RMSE. In the leading vehicle accelerating, the DRS has the lowest RMSE, while the RMSEs of the other two are relatively large. It is due to the speed risk surrogate model proposed in DRS. In the speed risk surrogate, the desired velocity of the object vehicle is introduced and affected by the velocity of traffic flow and interactive vehicles. In DRS, the simulated speed is better integrated into the surrounding vehicles and has less impact on the current traffic flow, which is in line with the driver's operating psychology and interaction behaviors.

Moreover, DRS’s RMSE of Trajectory No.1396 is 1.85, as shown in the top left corner of \mbox{Fig. \ref{1vha}}(c) and that of Trajectory No. 3044 is 3.04 in the top left of \mbox{Fig. \ref{2vha}}(c). The acceleration in No.1396 frequently changes, while that in No. 3044 changes less. That can be seen that the proposed method has better performance when the speed oscillates frequently.

\mbox{Fig. \ref{1vha}} (b) and \mbox{Fig. \ref{2vha}} (b) are acceleration curves obtained from the simulation of three models. The acceleration is derived from field strength and showed safe interactions more objectively. With the speed change of the vehicle ahead, the rear vehicle constantly changes its velocity to adapt to the front vehicle. It can be seen from \mbox{Fig. \ref{1vha}} (b) \uppercase\expandafter{\romannumeral2}, that the leading vehicle starts accelerating at 883.31s. DRS's acceleration becomes greater than zero at 883.74s. That is, DRS takes the lead in responding to the acceleration behavior. By comparison, IDM and DRPFM start accelerating at 885.15s and 883.94s, respectively, and the manual vehicle with real data starts accelerating at 885.51s. The three models are all able to accelerate in advance, and DRS is the earliest. The same comparison can be seen 
in \mbox{Fig. \ref{1vha}} (b) \uppercase\expandafter{\romannumeral4} and \uppercase\expandafter{\romannumeral6}, and \mbox{Fig. \ref{2vha}} (b) \uppercase\expandafter{\romannumeral2} where DRS and DRPFM start accelerating earlier than IDM in each scenario. Due to both DRS and DRPFM considering the acceleration of the front vehicle, they respond faster to the acceleration of the front vehicle than IDM. Moreover, DRS considers a speed risk surrogate, which makes it more sensitive to the front vehicle velocity than DRPFM during car-following interactions. 

When the leading vehicle starts decelerating, the rear vehicle needs to decelerate to avoid collisions. Let's take a look at the performance of the three algorithms when the front vehicle starts slowing down. It can be seen from \mbox{Fig. \ref{2vha}} (b) \uppercase\expandafter{\romannumeral3} that acceleration of the front vehicle is lower than zero at 887.67s, which means it starts decelerating. Similarly, it can be seen from the graph that the following manual-driven vehicle starts decelerating in 888.94s to avoid collision with the front vehicle. Compared with the acceleration curves of the three models, DRS and DRPFM perceive the risk and start slowing down at 888.01s, while IDM starts slowing down at 888.88s. The performance is consistent with that of the front vehicle when accelerating, which shows the response of DRS and DRPFM to the front vehicle deceleration is also more sensitive than IDM. In \mbox{Fig. \ref{1vha}} (b) \uppercase\expandafter{\romannumeral5}, it is earlier than DRPFM (899.18s) that DRS perceives the risk and starts slowing down at 898.71s. It can be seen that DRS has better adaptability and predictability in rear-end collision scenarios. 

\mbox{Fig. \ref{1vha}} (d) and \mbox{Fig. \ref{2vha}} (d) show relevant velocity curves. Responses to speed change also can be seen from another perspective. The sensitivity to other vehicle's motion changes brings benefits to the subject vehicle driving. An earlier response allows the following vehicle to take the initiative to ensure safety and efficiency, leaving more time for deceleration and acceleration adjustments. As \mbox{Fig. \ref{1vha}} (d) \uppercase\expandafter{\romannumeral6} and \mbox{Fig. \ref{2vha}} (d) \uppercase\expandafter{\romannumeral2} showing, benefiting from the sensitivity, the speed fitting effect of DRS is better than IDM and DRPFM in the car-following scenario when the front vehicle continues to accelerate. Interactive vehicle risk, restrictions risk, and speed risk are proposed and integrated into DRS. Instead of IDM and DRPFM considering velocity from one side, DRS study velocity influence more comprehensively in three aspects. Based on studying the front vehicle velocity in the vehicle energy, interactive vehicle risk surrogate, and speed risk surrogate, DRS shows the mechanism of the influence of the front vehicle speed on the risk perception and driving efficiency of the rear vehicle to a certain extent and makes a significant contribution to the research of the risk assessment model.

In addition, in traditional methods, IDM does not consider the impact of vehicle attributes. The virtual mass of the vehicle is considered in DRPFM. However, in the general naturalistic driving data, it is difficult to obtain the vehicle quality data. Therefore, the consideration of vehicle mass in models based on virtual mass is ideal and not enough to support practical applications. Although the driver in the rear vehicle cannot acquire the specific value of the quality of the front vehicle, he can intuitively gain the shape of the front vehicle. This paper select the naturalistic driving dataset from the aerial photography perspective. Vehicle sizes provide a realistic basis for the calibration of the virtual energy of DRS. 

Finally, it is worth discovering that the model proposed in this paper is calibrated by real data, while errors between simulation results and real trajectories still exist. But we found that DRS simulations perform better trajectory behavior than manual driving to some extent. In terms of space headway, when the vehicle in front starts decelerating, DRS usually maintains a larger space headway distance than manual driving vehicles to ensure more safety. In terms of velocity and acceleration, DRS always reacts to changes in the preceding vehicle earlier than human-driven vehicles. In addition, velocity fluctuations of DRS are generally smaller than that of manual driving, as shown in \mbox{Fig. \ref{1vha}} (d) and \mbox{Fig. \ref{2vha}} (d). That avoid unnecessary acceleration and deceleration and save fuel partly. Based on the above performances, DRS effectively captures the human perception of risk and the driving rules in the car-following case. Thus, it has the potential to provide solutions for planning driving trajectories and reduce the driving instability caused by human reaction time, driving operation level, and individual driving style.

\section{Conclusion}

In this study, DRS is proposed to reflect risks during safe interactions in driving. In comparison with existing driving risk potential models, the proposed model introduces virtual energy and speed risks into the model for the first time. The driving risk assessment issue is tackled from a new perspective that is an integrated model consisting of interactive vehicle risk, restrictions risk, and speed risk. DRS reflects the driving risk distribution of vehicles at different relative speeds and different following clearances and describes the general phenomenon more naturally. DRS model is applied to the car-following process. Through the calibration based on the naturalistic driving dataset, the car-following motion at the micro-level of the vehicle is simulated. The proposed model performs well in trajectory estimation, space headway estimation, and velocity estimation and even has better driving behaviors than the manual driving vehicle. Safe interaction in driving is always a core issue. The proposed DRS has the potential to be applied to guide the operation of appropriate countermeasures in intelligent transportation systems, connected and automated vehicle environments, or intelligent assistance systems.   

The academic significance of this paper includes two parts: one is the modeling of driving risk, and the other is the modeling of car-following. DRS proposed in this paper greatly improves the traditional risk potential field model. The designed virtual energy is more in line with the driver's intuitive feeling of interactive vehicles. Meanwhile, it also solves the problem that the vehicle quality cannot be obtained from naturalistic driving datasets of aerial photography, which leads to the idealized processing of virtual quality in conventional models. The virtual energy and speed risk surrogate considers the front vehicle speed and the influence mechanism of the front vehicle speed on risk perception. It presents the driving efficiency of the subject vehicle. This enables the proposed model to respond to the velocity variation earlier and provide a more comfortable and efficient trajectory service. In addition, based on the risk field strength, acceleration that can guide the vehicle's motion is outputted. On the one hand, acceleration is utilized to solve the problem that the threshold is not unified in different traffic scenarios. The risk surrogate model can be calibrated automatically based on acceleration according to the real driving dataset. This also simplifies the complexity of the calibration of the field strength threshold. On the other hand, the output acceleration realizes the modeling of car-following behavior and can be derived from the modeling of lane-changing behavior. A new foundation can be laid for the risk distribution in traffic flow through the connection between the microscopic model and the macroscopic model of traffic flow.

This study reveals a great potential for applying the proposed algorithm. However, this study has some limitations which need to be further improved. First, there is a lack of research on different driving styles. The parameters of the proposed model are calibrated based on the minimum error of the naturalistic driving dataset. It is expected that further research will be conducted on datasets with driver characteristics. The corresponding parameters of specific categories of drivers will be calibrated. This will enhance the anthropomorphic performance of the proposed model. Secondly, the traffic elements considered in this paper are relatively simple and have not been verified on urban roads. In further research, we will calibrate parameters in new scenarios and verify the car-following performance in more complex traffic scenarios. Nonetheless, we believe that the application of DRS and car-following in this paper may open a new perspective for many tasks to advance the development of micro-behavioral models and traffic flow models.




\bibliographystyle{alpha}
\bibliography{sample}

\newcommand{\etalchar}[1]{$^{#1}$}
\begin{thebibliography}{WAW{\etalchar{+}}15}

\bibitem[BK10]{blas2010equity}
Erik Blas and Anand~Sivasankara Kurup.
\newblock {\em Equity, social determinants and public health programmes}.
\newblock World Health Organization, 2010.

\bibitem[CWW{\etalchar{+}}17]{chen2017monitoring}
Faan Chen, Jianjun Wang, Jiaorong Wu, Xiaohong Chen, and P~Christopher Zegras.
\newblock Monitoring road safety development at regional level: A case study in
  the asean region.
\newblock {\em Accident Analysis \& Prevention}, 106:437--449, 2017.

\bibitem[FLWF21]{LiZhibin}
Ruyi Feng, Zhibin Li, Qifan Wu, and Changyan Fan.
\newblock Aerial video vehicle detection target correlation and space-time
  trajectory matching.
\newblock In {\em Traffic information and security(Chinese journal)}, pages
  39(2):61--69+77, 2021.

\bibitem[FS22]{fu_bayesian_2022}
Chuanyun Fu and Tarek Sayed.
\newblock Bayesian dynamic extreme value modeling for conflict-based real-time
  safety analysis.
\newblock {\em Analytic Methods in Accident Research}, 34:100204, 2022.

\bibitem[Gre93]{greenwade93}
George~D. Greenwade.
\newblock The {C}omprehensive {T}ex {A}rchive {N}etwork ({CTAN}).
\newblock {\em TUGBoat}, 14(3):342--351, 1993.

\bibitem[HEDJ13]{hojjati2013stochastic}
Khashayar Hojjati-Emami, Balbir~S Dhillon, and Kouroush Jenab.
\newblock The stochastic and integrative prediction methodology and modeling
  for reliability of pedestrian crossing on roads.
\newblock {\em Journal of Transportation Safety \& Security}, 5(3):257--272,
  2013.

\bibitem[HWF{\etalchar{+}}20]{huang2020probabilistic}
Heye Huang, Jianqiang Wang, Cong Fei, Xunjia Zheng, Yibin Yang, Jinxin Liu,
  Xiangbin Wu, and Qing Xu.
\newblock A probabilistic risk assessment framework considering lane-changing
  behavior interaction.
\newblock {\em Science China Information Sciences}, 63(9):1--15, 2020.

\bibitem[HWL12]{hsu2012conceptual}
Tien-Pen Hsu, Guo-Yu Weng, and Yu-Jui Lin.
\newblock Conceptual structure of a novel car-following model upon
  gravitational field concept.
\newblock In {\em 19th ITS World CongressERTICO-ITS EuropeEuropean
  CommissionITS AmericaITS Asia-Pacific}, 2012.

\bibitem[Kar10]{Karaboga:2010}
D.~Karaboga.
\newblock {A}rtificial bee colony algorithm.
\newblock {\em Scholarpedia}, 5(3):6915, 2010.
\newblock revision \ 91003.

\bibitem[KT86]{krogh1986integrated}
Bruce Krogh and Charles Thorpe.
\newblock Integrated path planning and dynamic steering control for autonomous
  vehicles.
\newblock In {\em Proceedings. 1986 IEEE International Conference on Robotics
  and Automation}, volume~3, pages 1664--1669. IEEE, 1986.

\bibitem[LGJ{\etalchar{+}}20]{li2020dynamic}
Linheng Li, Jing Gan, Xinkai Ji, Xu~Qu, and Bin Ran.
\newblock Dynamic driving risk potential field model under the connected and
  automated vehicles environment and its application in car-following modeling.
\newblock {\em IEEE Transactions on Intelligent Transportation Systems},
  23(1):122--141, 2020.

\bibitem[LGY{\etalchar{+}}20]{li2020risk}
Linheng Li, Jing Gan, Ziwei Yi, Xu~Qu, and Bin Ran.
\newblock Risk perception and the warning strategy based on safety potential
  field theory.
\newblock {\em Accident Analysis \& Prevention}, 148:105805, 2020.

\bibitem[LGZ{\etalchar{+}}20]{li2020novel}
Linheng Li, Jing Gan, Kun Zhou, Xu~Qu, and Bin Ran.
\newblock A novel lane-changing model of connected and automated vehicles:
  Using the safety potential field theory.
\newblock {\em Physica A: Statistical Mechanics and its Applications},
  559:125039, 2020.

\bibitem[LHvL{\etalchar{+}}21]{lu2021performance}
Chang Lu, Xiaolin He, Hans van Lint, Huizhao Tu, Riender Happee, and Meng Wang.
\newblock Performance evaluation of surrogate measures of safety with
  naturalistic driving data.
\newblock {\em Accident Analysis \& Prevention}, 162:106403, 2021.

\bibitem[LLZ{\etalchar{+}}20]{li2020threat}
Yang Li, Keqiang Li, Yang Zheng, Bernhard Morys, Shuyue Pan, and Jianqiang
  Wang.
\newblock Threat assessment techniques in intelligent vehicles: A comparative
  survey.
\newblock {\em IEEE Intelligent Transportation Systems Magazine}, 13(4):71--91,
  2020.

\bibitem[LM10]{lord2010statistical}
Dominique Lord and Fred Mannering.
\newblock The statistical analysis of crash-frequency data: A review and
  assessment of methodological alternatives.
\newblock {\em Transportation research part A: policy and practice},
  44(5):291--305, 2010.

\bibitem[LZLS18]{liu_multivariate_2018}
Chenhui Liu, Mo~Zhao, Wei Li, and Anuj Sharma.
\newblock Multivariate random parameters zero-inflated negative binomial
  regression for analyzing urban midblock crashes.
\newblock {\em Analytic Methods in Accident Research}, 17:32--46, 2018.

\bibitem[MB14]{mannering_analytic_2014}
Fred~L. Mannering and Chandra~R. Bhat.
\newblock Analytic methods in accident research: {Methodological} frontier and
  future directions.
\newblock {\em Analytic Methods in Accident Research}, 1:1--22, 2014.

\bibitem[MHZB21]{mohammadian_integrating_2021}
Saeed Mohammadian, Md~Mazharul Haque, Zuduo Zheng, and Ashish Bhaskar.
\newblock Integrating safety into the fundamental relations of freeway traffic
  flows: A conflict-based safety assessment framework.
\newblock {\em Analytic Methods in Accident Research}, 32:100187, 2021.

\bibitem[MRN13]{matsumi2013autonomous}
Ryosuke Matsumi, Pongsathorn Raksincharoensak, and Masao Nagai.
\newblock Autonomous braking control system for pedestrian collision avoidance
  by using potential field.
\newblock {\em IFAC Proceedings Volumes}, 46(21):328--334, 2013.

\bibitem[NBS16]{nadimi2016calibration}
Navid Nadimi, Hamid Behbahani, and Hamid~Reza Shahbazi.
\newblock Calibration and validation of a new time-based surrogate safety
  measure using fuzzy inference system.
\newblock {\em Journal of traffic and transportation engineering (English
  edition)}, 3(1):51--58, 2016.

\bibitem[Ni13]{ni2013unified}
Daiheng Ni.
\newblock A unified perspective on traffic flow theory part i: The field
  theory.
\newblock {\em Applied Mathematical Sciences}, 7(39):1929--1946, 2013.

\bibitem[NKI{\etalchar{+}}16]{norros2016palm}
Ilkka Norros, Pirkko Kuusela, Satu Innamaa, Eetu Pilli-Sihvola, and Riikka
  Rajam{\"a}ki.
\newblock The palm distribution of traffic conditions and its application to
  accident risk assessment.
\newblock {\em Analytic methods in accident research}, 12:48--65, 2016.

\bibitem[QKOJ14]{qu2014safety}
Xiaobo Qu, Yan Kuang, Erwin Oh, and Sheng Jin.
\newblock Safety evaluation for expressways: a comparative study for
  macroscopic and microscopic indicators.
\newblock {\em Traffic injury prevention}, 15(1):89--93, 2014.

\bibitem[RAMN13]{raksincharoensak2013predictive}
Pongsathorn Raksincharoensak, Yuta Akamatsu, Katsumi Moro, and Masao Nagai.
\newblock Predictive braking assistance system for intersection safety based on
  risk potential.
\newblock {\em IFAC Proceedings Volumes}, 46(21):335--340, 2013.

\bibitem[RHZ21]{sharifi2016time}
Rezaur Rahman, Samiul Hasan, and Mohamed~H. Zaki.
\newblock Towards reducing the number of crashes during hurricane evacuation:
  Assessing the potential safety impact of adaptive cruise control systems.
\newblock {\em Transportation Research Part C: Emerging Technologies},
  128:103188, 2021.

\bibitem[RKCL16]{rasekhipour2016potential}
Yadollah Rasekhipour, Amir Khajepour, Shih-Ken Chen, and Bakhtiar Litkouhi.
\newblock A potential field-based model predictive path-planning controller for
  autonomous road vehicles.
\newblock {\em IEEE Transactions on Intelligent Transportation Systems},
  18(5):1255--1267, 2016.

\bibitem[RS94]{reichardt1994collision}
Dirk Reichardt and Jeong Shick.
\newblock Collision avoidance in dynamic environments applied to autonomous
  vehicle guidance on the motorway.
\newblock In {\em Proceedings of the Intelligent Vehicles' 94 Symposium}, pages
  74--78. IEEE, 1994.

\bibitem[SB08]{sattel2008robotics}
Thomas Sattel and Thorsten Brandt.
\newblock From robotics to automotive: Lane-keeping and collision avoidance
  based on elastic bands.
\newblock {\em Vehicle System Dynamics}, 46(7):597--619, 2008.

\bibitem[SEE11]{sun2011elevated}
Carlos Sun, Praveen Edaram, and Kyle Ervin.
\newblock Elevated-risk work zone evaluation of temporary rumble strips.
\newblock {\em Journal of Transportation Safety \& Security}, 3(3):157--173,
  2011.

\bibitem[SS19]{schlogl2019methodological}
Matthias Schlögl and Rainer Stütz.
\newblock Methodological considerations with data uncertainty in road safety
  analysis.
\newblock {\em Accident Analysis \& Prevention}, 130:136--150, 2019.

\bibitem[SVG{\etalchar{+}}16]{staton2016road}
Catherine Staton, Joao Vissoci, Enying Gong, Nicole Toomey, Rebeccah Wafula,
  Jihad Abdelgadir, Yi~Zhou, Chen Liu, Fengdi Pei, Brittany Zick, et~al.
\newblock Road traffic injury prevention initiatives: a systematic review and
  metasummary of effectiveness in low and middle income countries.
\newblock {\em PloS one}, 11(1):e0144971, 2016.

\bibitem[TCW17]{tan2017use}
Dongkui Tan, Wuwei Chen, and Hongbo Wang.
\newblock On the use of monte-carlo simulation and deep fourier neural network
  in lane departure warning.
\newblock {\em IEEE Intelligent Transportation Systems Magazine}, 9(4):76--90,
  2017.

\bibitem[TVH01]{tsourveloudis2001autonomous}
Nikos~C Tsourveloudis, Kimon~P Valavanis, and Timothy Hebert.
\newblock Autonomous vehicle navigation utilizing electrostatic potential
  fields and fuzzy logic.
\newblock {\em IEEE transactions on robotics and automation}, 17(4):490--497,
  2001.

\bibitem[WAW{\etalchar{+}}15]{ward2015extending}
James~R Ward, Gabriel Agamennoni, Stewart Worrall, Asher Bender, and Eduardo
  Nebot.
\newblock Extending time to collision for probabilistic reasoning in general
  traffic scenarios.
\newblock {\em Transportation Research Part C: Emerging Technologies},
  51:66--82, 2015.

\bibitem[WHO20]{WHO}
WHO.
\newblock Road traffic injuries.
\newblock \url{https://www.who.int/health-topics/road-safety\#tab=tab_1}, 2020.

\bibitem[WLC10]{wang2010applying}
Jyun-Guo Wang, Cheng-Jian Lin, and Shyi-Ming Chen.
\newblock Applying fuzzy method to vision-based lane detection and departure
  warning system.
\newblock {\em Expert systems with applications}, 37(1):113--126, 2010.

\bibitem[WPH{\etalchar{+}}14]{paschalidis_deriving_2020}
Chaozhong Wu, Liqun Peng, Zhen Huang, Ming Zhong, and Duanfeng Chu.
\newblock A method of vehicle motion prediction and collision risk assessment
  with a simulated vehicular cyber physical system.
\newblock {\em Transportation Research Part C: Emerging Technologies},
  47:179--191, 2014.

\bibitem[WWL15]{wang2015driving}
Jianqiang Wang, Jian Wu, and Yang Li.
\newblock The driving safety field based on driver--vehicle--road interactions.
\newblock {\em IEEE Transactions on Intelligent Transportation Systems},
  16(4):2203--2214, 2015.

\bibitem[WWZ{\etalchar{+}}16]{wang2016driving}
Jianqiang Wang, Jian Wu, Xunjia Zheng, Daiheng Ni, and Keqiang Li.
\newblock Driving safety field theory modeling and its application in
  pre-collision warning system.
\newblock {\em Transportation research part C: emerging technologies},
  72:306--324, 2016.

\bibitem[WXX{\etalchar{+}}18]{alozi_evaluating_2022}
Chen Wang, Chengcheng Xu, Jingxin Xia, Zhendong Qian, and Linjun Lu.
\newblock A combined use of microscopic traffic simulation and extreme value
  methods for traffic safety evaluation.
\newblock {\em Transportation Research Part C: Emerging Technologies},
  90:281--291, 2018.

\bibitem[ZHGN18]{arun_systematic_2021}
Zhenhua Zhang, Qing He, Jing Gao, and Ming Ni.
\newblock A deep learning approach for detecting traffic accidents from social
  media data.
\newblock {\em Transportation Research Part C: Emerging Technologies},
  86:580--596, 2018.

\bibitem[ZHNX18]{zheng2018novel}
Xunjia Zheng, Bin Huang, Daiheng Ni, and Qing Xu.
\newblock A novel intelligent vehicle risk assessment method combined with
  multi-sensor fusion in dense traffic environment.
\newblock {\em Journal of intelligent and connected vehicles}, 2018.

\bibitem[ZHW{\etalchar{+}}21]{zheng2021behavioral}
Xunjia Zheng, Heye Huang, Jianqiang Wang, Xiaocong Zhao, and Qing Xu.
\newblock Behavioral decision-making model of the intelligent vehicle based on
  driving risk assessment.
\newblock {\em Computer-Aided Civil and Infrastructure Engineering},
  36(7):820--837, 2021.

\bibitem[ZZG{\etalchar{+}}18]{zheng2018novel2}
Xunjia Zheng, Di~Zhang, Hongbo Gao, Zhiguo Zhao, Heye Huang, and Jianqiang
  Wang.
\newblock A novel framework for road traffic risk assessment with hmm-based
  prediction model.
\newblock {\em Sensors}, 18(12):4313, 2018.

\end{thebibliography}

\end{document}